\documentclass{ifacconf}
\usepackage{amssymb}
\usepackage{caption}
\usepackage{natbib}    
\usepackage{graphicx}
\usepackage{epsfig}
\usepackage{subcaption}
\usepackage{multicol, blindtext}
\usepackage{lipsum}
\usepackage{multirow}
\usepackage{amsmath}
\usepackage{amsfonts}
\usepackage{amssymb}
\usepackage{bm}
\usepackage[usenames]{color}
\newtheorem{example}{Example}
\newtheorem{definition}{Definition}

\newtheorem{lemma}{Lemma}
\newtheorem{theorem}{Theorem}
\newtheorem{corollary}{Corollary}
\usepackage{algorithm,algorithmic}
\usepackage{float}
\floatstyle{boxed}
\restylefloat{figure}

\renewcommand{\vec}{\bm}
  
   \newtheorem{assumption}{Assumption}
\newtheorem{remark}{Remark}
\newcommand{\oomit}[1]{}

\begin{document}
\begin{frontmatter}

\title{Robust Regions of Attraction Generation for State-Constrained Perturbed Discrete-Time Polynomial Systems}

\author{Bai Xue and Naijun Zhan and Yangjia Li}

\address{University of Chinese Academy of Sciences, Beijing, China\\
State Key Lab. of Computer Science, Institute of Software, CAS, China\\ (e-mail:\{xuebai,znj,yangjia\}@ios.ac.cn)}

\begin{abstract}                        
In this paper we propose a convex programming based method for computing robust regions of attraction for state-constrained perturbed discrete-time polynomial systems. The robust region of attraction of interest is a set of states such that every possible trajectory initialized in it will approach an equilibrium state while never violating the specified state constraint, regardless of the actual perturbation. Based on a Bellman equation which characterizes the interior of the maximal robust region of attraction as the strict one sub-level set of its unique bounded and continuous solution, we construct a semi-definite program for computing robust regions of attraction.  Under appropriate assumptions, the existence of solutions to the constructed semi-definite program is guaranteed and there exists a sequence of solutions such that their strict one sub-level sets inner-approximate and converge to  the interior of the maximal robust region of attraction in measure. Finally, we demonstrate the method by two examples.
\end{abstract}


\end{frontmatter}

\section{Introduction}
\label{Int}
Discrete-times systems, which are governed by difference equations or iterative processes, may result from discretizing continuous systems or modeling evolution systems for which the time scale is discrete. They are prevalent in signal processing, population dynamics, scientific computation and so forth, e.g., \citep{kot1986}. The polynomial discrete-time systems are the type of systems whose dynamics are described in polynomial forms. This system is classified as an important class of nonlinear systems due to the fact that many nonlinear systems can be modelled as, transformed into, or approximated by polynomial systems, e.g., \citep{halanay2000stability}.   

A fundamental problem in control engineering consists of determining the robust region of attraction of an equilibrium \citep{slotine1991}, which is a set of states such that every trajectory starting from it will approach this equilibrium while never leaving a specified state constraint set irrespective of the actual perturbation. Its applications include biology systems \citep{merola2008} and ecology systems \citep{ludwig1997}, among others. Computing robust regions of attraction has been the subject of extensive research over the past several decades, resulting in the emergence of many computational approaches, e.g., Lyapunov function based methods \citep{zubov1964,lasaal1961,coutinho2013,genesio1985,giesl2014},  trajectory reversing methods \citep{genesio1985} and moment-based methods \citep{henrion2013convex,korda2013inner}. 

Lyapunov function based methods are still dominant in estimating robust regions of attraction \citep{khalil2002nonlinear}. Generally, the search for Lyapunov functions is non-trivial for nonlinear systems due to the
non-constructive nature of the Lyapunov theory, apart from some cases where the Jacobian matrix of the linearized system associated with the nonlinear system of interest
is Hurwitz. However, with the advance of real algebraic geometry  and polynomial optimization in the last decades, especially the sum-of-squares (SOS) decomposition technique \citep{parrilo2000}, finding a Lyapunov function which is decreasing over a given state constraint set could be reduced to a convex programming problem for polynomial
systems \citep{papachristodoulou2002}. This results in a large amount of findings which adopt convex optimization based approaches to the search for polynomial Lyapunov functions, e.g., \citep{anderson2015}. However, if we return to the problem of estimating robust
domains of attraction, it resorts to addressing a bilinear semi-definite program, e.g., \citep{jarvis2003,tan2008}, which falls within the non-convex programming framework and is notoriously hard to solve. Also, the existence of polynomial solutions to (bilinear) semi-definite
programs is not explored  in the literature, especially for perturbed systems.
 
In this paper we propose a novel semi-definite programming based method for computing robust regions of attraction for state-constrained perturbed discrete-time polynomial systems with an equilibrium state, which is uniformly locally exponentially stable. It is worth remarking here that the method proposed in this paper can also be applied to the computation of robust regions of attraction for polynomial systems with an asymptotically stable equilibrium state, as highlighted in Remark \ref{asym}. The semi-definite program is constructed by relaxing a modified Bellman equation which characterizes the interior of the maximal robust region of attraction as the strict one sub-level set of its unique bounded and continuous solution \citep{xue2018regions}. It falls within the convex programming framework and can be solved efficiently in polynomial time via interior-point methods. Moreover, the existence of solutions to the constructed semi-definite program is guaranteed and there exists a sequence of polynomial solutions such that their strict one sub-level sets inner-approximate and converge to the interior of the maximal robust region of attraction  in measure under appropriate assumptions. Finally, we demonstrate our method by two examples. 

The closely related works to the present one in spirit are \citep{summers2013,henrion2013convex,korda2013inner, xue2018under,xue2019inner,WZF19}. The work in \citep{summers2013} employed semi-definite programs to solve discrete-time stochastic optimal control problems by relaxing the Bellman equation. The works in \citep{henrion2013convex,korda2013inner} respectively considered outer and inner approximations of the maximal region of attraction over finite time horizons. However, the present one focuses on inner-approximations of the maximal region of attraction over the infinite time horizon. Recently, semi-definite programming based methods were proposed in \citep{xue2018under,xue2019inner,WZF19} for computing reachable sets of continuous-time polynomial systems by relaxing Hamilton-Jacobi equations, where the computation of reachable sets over finite time horizons was studied in \citep{xue2018under,xue2019inner} and the computation of robust invariant sets over the infinite time horizon was studied in \citep{WZF19}. Trajectories starting from the robust invariant set in \citep{WZF19} are not required to approach an equilibrium. Also, the existence of solutions to the constructed semi-definite program in \citep{WZF19} is not guaranteed. In contrast, the present work considers the computation of robust regions of attraction over the infinite time horizon for discrete-time polynomial systems by relaxing Bellman equations. Trajectories starting from the robust region of attraction are required to approach an equilibrium. Moreover, the existence of solutions to the constructed semi-definite program in the present work is guaranteed under appropriate conditions.

This paper is structured as follows. In Section \ref{Pre} basic notions and the problem of interest are introduced. After presenting our method for synthesizing robust regions of attraction in Section \ref{IG}, we evaluate it on two examples in Section \ref{ex}. Finally, we conclude this paper in Section \ref{con}.

\section{Preliminaries}\label{Pre}
In this section we describe the system of interest and the concept of robust regions of attraction. 

The notions will be used in this paper: $\mathbb{R}^n$ denotes the set of $n$-dimensional real vectors. $\mathbb{R}[\cdot]$ denotes the ring of polynomials with real coefficients in variables given by the argument. $\mathbb{R}_k[\cdot]$ denotes the vector space of real multivariate polynomials of total degree $\leq k$. $\Delta^{\circ}$, $\partial \Delta$, $\overline{\Delta}$ and $\Delta^{c}$ denote the interior, boundary, closure and complement of a set $\Delta$, respectively. The difference of two sets $A$ and $B$ is denoted by $A\setminus B$.  $\mu(A)$ denotes the Lebesgue measure on $A\subset \mathbb{R}^n$. $\mathbb{N}$ denotes the set of non-negative integers. $\|\bm{x}\|$ denotes the 2-norm, i.e., $\|\bm{x}\|=\sqrt{\sum_{i=1}^n x_i^2}$, where $\bm{x}=(x_1,\ldots,x_n)^{\top}$. $B(\bm{0},r)$ denotes a ball of radius $r>0$ and center $\bm{0}$, i.e., $B(\bm{0},r)=\{\bm{x}\mid \|\bm{x}\|^2\leq r\}$. Vectors are denoted by boldface letters.

The perturbed discrete-time system of interest in this paper is of the following form
\begin{equation}
\label{systems}
\bm{x}(k+1)=\bm{f}(\bm{x}(k),\bm{d}(k)), k\in \mathbb{N},
\end{equation}
where $\bm{x}(\cdot):\mathbb{N}\rightarrow \mathbb{R}^n$, $\bm{d}(\cdot):\mathbb{N}\rightarrow D$, $$D=\{\bm{d}\in \mathbb{R}^m\mid \wedge_{i=1}^{m_d} [h_i^{D}(\bm{d})\leq 0]\}$$ is a compact semi-algebraic subset in $\mathbb{R}^m$ with $h_i^{D}\in \mathbb{R}[\bm{d}]$, and $\bm{f}\in \mathbb{R}[\bm{x},\bm{d}]$ with $\bm{f}(\bm{0},\bm{d})=\bm{0}$ for $\bm{d}\in D$.

In order to define our problem succinctly, we present the definition of a perturbation input policy $\pi$. 
\begin{definition}
\label{policy}
A perturbation input policy, denoted by $\pi$, refers to a function $\pi(k):\mathbb{N}\rightarrow D$. In addition, we denote the set of all perturbation input policies by $\mathcal{D}$.
\end{definition}

Given a perturbation input policy $\pi$, a trajectory of system \eqref{systems} is presented in Definition \ref{tra}.
\begin{definition}
\label{tra}
Given a perturbation input policy $\pi\in \mathcal{D}$, a trajectory of system \eqref{systems} initialized in $\bm{x}_0 \in \mathbb{R}^n$ is defined as $\bm{\phi}_{\bm{x}_0}^{\pi}(\cdot): \mathbb{N} \rightarrow \mathbb{R}^n$, where $\bm{\phi}_{\bm{x}_0}^{\pi}(0)=\bm{x}_0$ and
\begin{equation}
\begin{split}
&\bm{\phi}_{\bm{x}_0}^{\pi}(k+1)=\bm{f}(\bm{\phi}_{\bm{x}_0}^{\pi}(k),\pi(k)), \forall k\in \mathbb{N}.\\
\end{split}
\end{equation}
\end{definition}

We assume that $\bm{0}$ is uniformly locally exponentially stable. 
\begin{assumption}
\label{app}
The equilibrium state $\bm{0}$ is uniformly locally exponentially stable for system \eqref{systems}, i.e., there exist positive constants $M>0$, $r>0$ and $0<\lambda<1$ such that
 \begin{equation*}
 \label{r}
\begin{split}
\|\bm{\phi}_{\bm{x}_0}^{\pi}(k)\|\leq \lambda^k M\|\bm{x}_0\|,\forall \bm{x}_0\in B(\bm{0},r), \forall \pi\in \mathcal{D},\forall k\in \mathbb{N},
 \end{split}
 \end{equation*}
  where $B(\bm{0},r)\subset X$ and $X\subset\mathbb{R}^n$ is a state constraint set, which will be defined later.
\end{assumption}

Assumption \ref{app} implies the existence of a positive constant $\overline{\epsilon}$ such that $B(\bm{0},\overline{\epsilon})\subseteq X$ and
 \begin{equation}
 \label{half}
\bm{\phi}_{\bm{x}_0}^{\pi}(k)\in B(\bm{0},\frac{r}{2}), \forall \bm{x}_0\in B(\bm{0},\overline{\epsilon}), \forall k\in \mathbb{N}, \forall \pi\in \mathcal{D}.
\end{equation}
Since $0<\lambda<1$ in Assumption \ref{app}, $\overline{\epsilon}$ in \eqref{half} exists and can take the value of $\min\{\frac{r}{2M},\frac{r}{2}\}$.

Suppose that the state constraint set $$X=\{\bm{x}\in \mathbb{R}^n\mid \wedge_{i=1}^{n_X}[h^X_i(\bm{x})<1]\}$$ is a bounded and open set with $h_i^X(\bm{x})\in \mathbb{R}[\bm{x}]$. Also, $h_i^X(\bm{x})>0$ for $\bm{x}\neq \bm{0}$ and $h_i^X(\bm{0})=0$, $i=1,\ldots,n_X$. We formally define robust regions of attraction.
\begin{definition}[Robust Regions of Attraction]
\label{prob}
The maximal robust region of attraction $\mathcal{R}$ is the set of states such that every possible trajectory of system \eqref{systems} starting from it will approach the equilibrium state $\bm{0}$ while never leaving the state constraint set $X$, i.e.
\begin{equation}
\mathcal{R}=\left\{\bm{x}_0 \middle|\;
\begin{aligned}
&\bm{\phi}_{\bm{x}_0}^{\pi}(k)\in X, \forall k\in \mathbb{N}, \forall \pi\in \mathcal{D}, \\
&~\text{ and }\lim_{k\rightarrow \infty} \bm{\phi}_{\bm{x}_0}^{\pi}(k)=\bm{0}, \forall \pi\in \mathcal{D}
\end{aligned}
\right\}.
\end{equation}
Correspondingly, a robust region of attraction is a subset of the maximal robust region of attraction $\mathcal{R}$.
\end{definition}

\section{Robust Regions of Attraction Generation}
\label{IG}
In this section we present our semi-definite programming based method for computing robust regions of attraction by relaxing Bellman equations. Furthermore, we show that there exists a sequence of solutions to the semi-definite program such that their strict one sub-level sets can inner-approximate the interior of the maximal robust region of attraction in measure under appropriate assumptions. 

\subsection{Bellman Equations}
In this subsection we introduce a modified Bellman equation, to which the strict one sub-level set of the unique bounded and continuous solution is equal to the interior of the maximal robust region of attraction.

\begin{theorem} 
\label{bell}
The interior of the maximal robust region of attraction $\mathcal{R}$ is equal to the strict one sub-level set of the unique bounded and continuous solution $v(\bm{x}):\mathbb{R}^n \rightarrow [0,1]$ to the Bellman equation
\begin{equation}
\label{e-v}
\begin{cases}
&\min\big\{\inf_{\bm{d}\in D}\{v-v(\bm{f})-g\cdot(1-v)\},\\
&v-1+\min_{j\in \{1,\ldots,n_{X}\}} l(1-h^X_j)\big\}=0, \forall \bm{x}\in \mathbb{R}^n,\\
&v(\bm{0})=0,\\
\end{cases}
\end{equation}
where $g(\cdot):\mathbb{R}^n\rightarrow \mathbb{R}$ is a non-negative polynomial satisfying that $g(\bm{x})=0$ iff $\bm{x}=\bm{0}$, and 
$l(\cdot):\mathbb{R}\rightarrow \mathbb{R}$ with 
 \begin{equation}
\begin{split}
    l(x)=
    \left\{
                \begin{array}{lll}
                  x,& \text{if } x\geq 0, \\
                  0, &\text{otherwise.}\\                  
                \end{array}
              \right.
              \end{split}
\end{equation}
That is, $\mathcal{R}^{\circ}=\{\bm{x}\in \mathbb{R}^n \mid v(\bm{x})<1\}$.
\end{theorem}

The Bellman equation \eqref{e-v}  in Theorem \ref{bell} is a discrete-time version of Zubov's equation for state-constrained continuous-time systems in \citep{grune2015}, and can be constructed  by following the
reasoning in \citep{grune2015}. Its detailed derivation is shown in Appendix.

A direct consequence of Theorem \ref{bell} is that  if a continuous function $u(\cdot): \mathbb{R}^n\rightarrow \mathbb{R}$ satisfies \eqref{e-v}, then $u(\bm{x})$ satisfies the constraints: 
\begin{equation}
\label{upper1}
\left\{
\begin{array}{lll}
u-u(\bm{f})-g\cdot(1-u)\geq 0, &\forall \bm{x}\in \mathbb{R}^n, \forall \bm{d}\in D,\\
u-1+l(1-h^X_j)\geq 0, &\forall \bm{x}\in \mathbb{R}^n,\\ 
j=1,\ldots,n_{X}.\\
\end{array}
\right.
\end{equation}

\begin{corollary}
\label{upper}
Suppose a continuous function $u(\bm{x}):\mathbb{R}^n\rightarrow \mathbb{R}$ is a solution to \eqref{upper1}, then $u(\bm{x})\geq v(\bm{x})$, where $v(\bm{x})$ is the unique bounded and continuous solution to the Bellman equation \eqref{e-v}. Consequently, $\{\bm{x}\in \mathbb{R}^n\mid u(\bm{x})<1\}\subset \mathcal{R}^{\circ}$ and thus
$\{\bm{x}\in \mathbb{R}^n\mid u(\bm{x})<1\}$ is a robust region of attraction. 
\end{corollary}
\begin{pf}
The second constraint in \eqref{upper1} implies that $u(\bm{x})\geq 0$ for $\bm{x}\in \mathbb{R}^n$. 

Assume that there exists $\bm{y}_0\in \mathbb{R}^n$ such that $u(\bm{y}_0)<v(\bm{y}_0)$. 

 First let's assume $v(\bm{y}_0)\geq 1$. Obviously, $\bm{y}_0\neq \bm{0}$ and consequently $g(\bm{y}_0)>0$. Since $u$ satisfies \eqref{upper1} and $v(\bm{y}_0)>u(\bm{y}_0)$,  we have that $v(\bm{y}_0)-1+\min_{j\in \{1,\ldots,n_X\}}l(1-h^X_j(\bm{y}_0))>0$. Also, since $v$ satisfies \eqref{e-v}, we have that $$\inf_{\bm{d}\in D}\{v(\bm{y}_0)-v(\bm{f}(\bm{y}_0,\bm{d}))-g(\bm{y}_0)(1-v(\bm{y}_0))\}=0.$$ Since $v$ is continuous over $\mathbb{R}^n$ and $\bm{f}$ is continuous over $\mathbb{R}^n\times D$, there exists $\bm{d}'_1\in D$ such that $v(\bm{y}_0)-v(\bm{f}(\bm{y}_0,\bm{d}'_1))-g(\bm{y}_0)(1-v(\bm{y}_0))=0$. Since $u(\bm{y}_0)-u(f(\bm{y}_0,\bm{d}'_1))-g(\bm{y}_0)(1-u(\bm{y}_0))\geq 0$, we obtain that 
\begin{equation*}
\begin{split}
v(\bm{f}(\bm{y}_0,\bm{d}'_1))-&u(\bm{f}(\bm{y}_0,\bm{d}'_1))\geq (v(\bm{y}_0)-u(\bm{y}_0))(1+g(\bm{y}_0)).
\end{split}
\end{equation*}
Let $\bm{y}_1=\bm{\phi}_{\bm{y}_0}^{\pi_{1}}(1)$, where $\pi_1(0)=\bm{d}'_1$, then $v(\bm{y}_1)>u(\bm{y}_1)$. Also, we have $v(\bm{y}_0)\leq v(\bm{y}_1)$. Moreover, $\bm{y}_1 \neq \bm{0}$, $g(\bm{y}_1)>0$. We continue the above deduction from $\bm{y}_0$ to $\bm{y}_1$, and obtain that there exists $\bm{d}'_2 \in D$ such that 
\begin{equation*}
\begin{split}
v(\bm{f}(\bm{y}_1,\bm{d}'_2))-&u(\bm{f}(\bm{y}_1,\bm{d}'_2))\geq (v(\bm{y}_1)-u(\bm{y}_1))(1+g(\bm{y}_1)).
\end{split}
\end{equation*}
Thus, we have 
\begin{equation*}
\begin{split}
&v(\bm{f}(\bm{y}_1,\bm{d}'_2))-u(\bm{f}(\bm{y}_1,\bm{d}'_2))\geq \\
&(v(\bm{y}_0)-u(\bm{y}_0))(1+g(\bm{y}_1))(1+g(\bm{y}_0)).
\end{split}
\end{equation*}
Let $\bm{y}_2=\bm{\phi}_{\bm{y}_1}^{\pi_{2}}(1)$, where $\pi_2(0)=\bm{d}'_2$, then $v(\bm{y}_2)>u(\bm{y}_2)$. Also, $v(\bm{y}_1)\leq v(\bm{y}_2)$.

Analogously, we deduce that for $k\in \mathbb{N}$,
\begin{equation*}
\begin{split}
v(\bm{f}&(\bm{y}_k,\bm{d}'_{k+1}))-u(\bm{f}(\bm{y}_k,\bm{d}'_{k+1}))\geq\\
& (v(\bm{y}_0)-u(\bm{y}_0))(1+g(\bm{y}_k))\cdots(1+g(\bm{y}_0)).
\end{split}
\end{equation*}
Moreover, let $\bm{y}_{k+1}=\bm{\phi}_{\bm{y}_k}^{\pi_{{k+1}}}(1)$, then $v(\bm{y}_k)\leq v(\bm{y}_{k+1})$, where $\pi_{k+1}(0)=\bm{d}'_{k+1}$. This implies that $\lim_{k\rightarrow \infty} \bm{y}_k\neq \bm{0}$ and thus $\bm{y}_k\notin B(\bm{0},\overline{\epsilon})$ for $k\in \mathbb{N}$, where $B(\bm{0},\overline{\epsilon})$ is defined in \eqref{half}.
Assume  $c_0=\inf\{g(\bm{x})\mid \bm{x}\in\mathbb{R}^n\setminus B(\bm{0},\overline{\epsilon})\}$. Clearly, $c_0>0$.  Such $c_0$ exists since $g(\bm{x})$ is a non-negative polynomial over $\mathbb{R}^n$ and $g(\bm{x})=0$ iff $\bm{x}=\bm{0}$. Therefore, 
\begin{equation*}
\begin{split}
v(\bm{f}(\bm{y}_k,\bm{d}'_{k+1}))-&u(\bm{f}(\bm{y}_k,\bm{d}'_{k+1}))\\
&\geq (v(\bm{y}_0)-u(\bm{y}_0))(1+c_0)^{k+1},
\end{split}
\end{equation*} implying that $\lim_{k\rightarrow \infty}v(\bm{y}_k)=\infty$, which contradicts the fact that $v$ is bounded over $\mathbb{R}^n$.  Thus, $v(\bm{y}_0)\leq u(\bm{y}_0)$.

Next, assume $v(\bm{y}_0)>u(\bm{y}_0)$ and $v(\bm{y}_0)<1$. According to Theorem \ref{bell}, every possible trajectory starting from $\bm{y}_0$ will eventually approach $\bm{0}$. Also, we have $$\inf_{\bm{d}\in D}\{v(\bm{y}_0,\bm{d})-v(\bm{f}(\bm{y}_0,\bm{d}))-g(\bm{y}_0)(1-v(\bm{y}_0))\}=0.$$ Following the deduction mentioned above, we have $$v(\bm{y}_k)-u(\bm{y}_{k})\geq v(\bm{y}_0)-u(\bm{y}_0),\forall k\in \mathbb{N}.$$ Therefore, we have that $\lim_{k\rightarrow \infty}v(\bm{y}_k)\geq v(\bm{y}_0)-u(\bm{y}_0)$, contradicting $\lim_{k\rightarrow \infty}v(\bm{y}_k)=0$. Thus, $v(\bm{y}_0) \leq u(\bm{y}_0)$. 

Therefore, $v(\bm{x})\leq u(\bm{x})$ for $\bm{x}\in \mathbb{R}^n$. Also,  since $\mathcal{R}^{\circ}=\{\bm{x}\in \mathbb{R}^n \mid v(\bm{x})<1\}$,  $\{\bm{x}\in \mathbb{R}^n \mid u(\bm{x})<1\}\subset \mathcal{R}^{\circ}$ holds.\qed
\end{pf}

From Corollary \ref{upper} we observe that a robust region of attraction can be found by solving \eqref{upper1} instead of \eqref{e-v}. 

\subsection{Semi-definite Programming Relaxation}
\label{sdpr}
In this subsection we construct a semi-definite program to compute robust regions of attraction based on \eqref{upper1}.  We observe that $u(\bm{x})$ is required to satisfy \eqref{upper1} over $\mathbb{R}^n$, which is a strong condition. Regarding this issue, we only consider \eqref{upper1} on the set $B(\bm{0},R)$, where the set $B(\bm{0},R)$ is defined in Assumption \ref{assump2}. In addition, we introduce another set $X_{\infty}$ with $X_{\infty}\subset X$ from Assumption \ref{assump2}, which is also defined in Assumption \ref{assump2}.  

 \begin{assumption}
 \label{assump2}
 \begin{enumerate}
 \item[(a)] $B(\bm{0},R)=\{\bm{x}\in \mathbb{R}^n\mid h_0(\bm{x})\leq R\}$, where $h_0(\bm{x})=\sum_{i=1}^nx_{i}^2$ and $R$ is a positive constant such that $\Omega(X)\subset B(\bm{0},R)$ with $\Omega(X)$ being the set of states reachable from the set $X$ within one step for system \eqref{systems}, i.e., 
$\Omega(X)=\{\bm{x}\in \mathbb{R}^n\mid \bm{x}=\bm{f}(\bm{x}_0,\bm{d}),\bm{x}_0\in X,\bm{d}\in D\}\cup X$.
 \item[(b)]  the set $X_{\infty}=\{\bm{x}\in \mathbb{R}^n\mid h_{\infty}(\bm{x})< 1\}$ is a robust region of attraction, where $h_{\infty}\in \mathbb{R}[\bm{x}]$. Besides, we assume that  $\bm{0}\in X_{\infty}^{\circ}$. It could be regarded as an initial conservative estimate of the maximal robust region of attraction.
\end{enumerate}
\end{assumption}
The set $X_{\infty}$ satisfies Assumption \ref{assump2} if $h_{\infty}$ is a (local) Lyapunov function for system \eqref{systems}. There are various methods for computing $h_{\infty}$, e.g., semi-definite programming based methods \citep{giesl2015review} and linear programming based methods \citep{giesl2014}. Also, the set $B(\bm{0},R)$ can be computed by solving a semi-definite programming problem as in \citep{magron2017}. In this paper, we assume both $X_{\infty}$ and $B(\bm{0},R)$ are given and thus their computations are not the focus of this paper.

Based on the sets $B(\bm{0},R)$ and $X_{\infty}$ in Assumption \ref{assump2}, we further relax constraint \eqref{upper1} and restrict the search for a continuous function $u(\bm{x})$ in the compact set $B(\bm{0},R)$, resulting in the following constraints:
\begin{equation}
\label{constraint0}
\begin{split}
&u-u(\bm{f})-g\cdot(1-u)\geq 0, \forall \bm{x}\in \overline{B(\bm{0},R)\setminus X_{\infty}},\forall \bm{d}\in D,\\
&u-1\geq 0, \forall \bm{x}\in B(\bm{0},R)\setminus X,\\
&u-h^X_j\geq 0, \forall \bm{x}\in \overline{X},\\
&j=1,\ldots,n_{X}. \\
\end{split}
\end{equation}
Obviously, $v(\bm{x})$ in \eqref{e-v} satisfies \eqref{constraint0}.

When the solution to \eqref{constraint0} is restricted to a polynomial, based on the sum-of-squares decomposition for multivariate polynomials, \eqref{upper1} could be reduced as the following sum-of-squares program \eqref{sos}.
\begin{algorithm}
\begin{equation}
\label{sos}
\begin{split}
&p_k^*=\inf \bm{w}\cdot \bm{l}\\
&{\rm s.t.}\\
&u_k-u_k(\bm{f})-g\cdot(1-u_k)=s_0+s_1\cdot(R-h_0)\\
&~~~~~~~~~~~~~~~~~~~~~~~~~~+s_2\cdot(h_{\infty}-1)-\sum_{i=1}^{m_d}s_{3,i}\cdot h_i^D,\\
&u_k-1=s_{4,j}+s_{5,j}\cdot(R-h_0)+s_{6,j}\cdot(h^X_j-1),\\
&u_k-h^X_j=s_{7,j}+s_{8,j}\cdot(R-h_0)+\sum_{l=1}^{n_{X}}s_{9,l,j}\cdot(1-h^X_l),\\
&j=1,\ldots,n_{X},\\
\end{split}
\end{equation}
where $\bm{w}\cdot \bm{l}=\int_{B(\bm{0},R)}u_k(\bm{x})d\bm{x}-\int_{X_{\infty}}u_k(\bm{x})d\bm{x}$, $\bm{l}$ is the vector of the moments of the Lebesgue measure over $\overline{B(\bm{0},R)\setminus X_{\infty}}$ indexed in the same basis in which the polynomial $u_k$ with coefficients $\bm{w}$ is expressed. The minimum is over the polynomial $u_k(\bm{x})\in R_k[\bm{x}]$ and the sum-of-squares polynomials $s_i(\bm{x},\bm{d})$, $i=0,\ldots,2$, $s_{3,i}(\bm{x},\bm{d})$, $i=1,\ldots,m_d$, $s_{i,j}(\bm{x})$, $s_{9,l,j}(\bm{x})$, $i=4,\ldots,8$, $j,l=1,\ldots,n_{X}$. 
\label{alg2}
\end{algorithm}

\begin{theorem}
\label{inner}
Under Assumption \ref{assump2}, if $u(\bm{x})\in \mathbb{R}_k[\bm{x}]$ is a solution to \eqref{sos}, then $\{\bm{x}\in B(\bm{0},R)\mid u(\bm{x})<1\}$ is a robust region of attraction.
\end{theorem}
\begin{pf}
According to the second constraint in \eqref{sos}, we
have $u(\bm{x})\geq 1$ for $\bm{x}\in B(\bm{0},R)\setminus X$. Therefore, $\{\bm{x}\in B(\bm{0},R)\mid u(\bm{x})<1\}\subset X$. Next we prove that every possible trajectory initialized in the set $\{\bm{x}\in B(\bm{0},R)\mid u(\bm{x})<1\}$ will approach the equilibrium state $\bm{0}$ eventually while never leaving the state constraint set $X$.

Assume that there exists $\bm{y}_0\in \{\bm{x}\in B(\bm{0},R)\mid u(\bm{x})<1\}$ and a perturbation input policy $\pi'$ such that $\bm{\phi}_{\bm{y}_0}^{\pi'}(k) \in X$ for $k=0,\ldots,k_0$ and $\bm{\phi}_{\bm{y}_0}^{\pi'}(k_0+1)\notin X$. It is obvious that $\bm{\phi}_{\bm{y}_0}^{\pi'}(k)\in X\setminus X_{\infty}$ for $k=0,\ldots,k_0$ since $X_{\infty}$ is a robust region of attraction. Since $\Omega(X)\subseteq B(\bm{0},R)$, where $\Omega(X)$ is defined in Assumption \ref{assump2}, $\bm{\phi}_{\bm{y}_0}^{\pi'}(k_0+1)\in B(\bm{0},R)\setminus X$, thus we obtain that 
\begin{equation}
\label{contra1}
u(\bm{\phi}_{\bm{y}_0}^{\pi'}(k_0+1))\geq 1.
\end{equation} 
However, since $\bm{\phi}_{\bm{y}_0}^{\pi'}(k)\in B(\bm{0},R)\setminus X_{\infty}$ for $k=0,\ldots, k_0+1$ and $u(\bm{y}_0)<1$, from the first constraint in \eqref{sos}, we have $$u(\bm{\phi}_{\bm{y}_0}^{\pi'}(k_0+1))<1,$$ contradicting \eqref{contra1}. Thus, every possible trajectory initialized in $\{\bm{x}_0\in B(\bm{0},R)\mid u(\bm{x}_0)<1\}$ never leaves $X$.

Lastly, we prove that every possible trajectory initialized in $\{\bm{x}\in B(\bm{0},R)\mid u(\bm{x})<1\}$ will  approach the equilibrium state $\bm{0}$ eventually. Since every possible trajectory initialized in the set $X_{\infty}$ will approach the equilibrium state $\bm{0}$ eventually, it is enough to prove that every possible trajectory initialized in the set $\{\bm{x}\in B(\bm{0},R)\mid u(\bm{x})<1\}\setminus X_{\infty}$ will enter the set $X_{\infty}$ in finite time. Assume that there exist $\bm{y}_0\in \{\bm{x}\in B(\bm{0},R)\mid u(\bm{x})<1\}$ and a perturbation input policy $\pi'$ such that $\bm{\phi}_{\bm{y}_0}^{\pi'}(k) \notin X_{\infty}, \forall k\in \mathbb{N}.$ Since $\bm{\phi}_{\bm{y}_0}^{\pi'}(k)\in X$ for $k\in \mathbb{N}$ and $u(\bm{x})\geq 0$ for $\bm{x}\in X$ (The fact that $u(\bm{x})\geq 0$ for $\bm{x}\in X$ can be obtained from the third constraint in \eqref{sos}.), $$u(\bm{\phi}_{\bm{y}_0}^{\pi'}(k))\geq 0, \forall k\in \mathbb{N}.$$ Moreover, $u(\bm{\phi}_{\bm{y}_0}^{\pi'}(k))<1$ holds for $k\in \mathbb{N}$. According to the first constraint in \eqref{sos}, we have $$u(\bm{\phi}_{\bm{y}_0}^{\pi'}(k))-u(\bm{\phi}_{\bm{y}_0}^{\pi'}(k+1))\geq g(\bm{\phi}_{\bm{y}_0}^{\pi'}(k))\cdot (1-u(\bm{\phi}_{\bm{y}_0}^{\pi'}(k)))$$ for $k\in \mathbb{N}$. Therefore, $$u(\bm{\phi}_{\bm{y}_0}^{\pi'}(k+1))\leq u(\bm{\phi}_{\bm{y}_0}^{\pi'}(k))-g(\bm{\phi}_{\bm{y}_0}^{\pi'}(k))(1-u(\bm{\phi}_{\bm{y}_0}^{\pi'}(k)))$$ and thus $$u(\bm{\phi}_{\bm{y}_0}^{\pi'}(k))\geq u(\bm{\phi}_{\bm{y}_0}^{\pi'}(k+1))$$ for $k\in \mathbb{N}$. Since $g(\bm{x})\in \mathbb{R}[\bm{x}]$ is positive over $\bm{x}\neq \bm{0}$, we obtain that $g(\bm{x})$ can attain a minimum over the compact set $\overline{X\setminus X_{\infty}}$. Let $$\epsilon'=\min_{\bm{x}\in \overline{X\setminus X_{\infty}}} g(\bm{x}),$$ it is obvious that $\epsilon'>0$. Therefore, we have $$u(\bm{\phi}_{\bm{y}_0}^{\pi'}(k+1))\leq u(\bm{\phi}_{\bm{y}_0}^{\pi'}(k))-\epsilon'(1-u(\bm{y}_0)), \forall k\in \mathbb{N}.$$ Therefore, $$u(\bm{\phi}_{\bm{y}_0}^{\pi'}(k+1))\leq u(\bm{y}_0)-(k+1)\epsilon'(1-u(\bm{y}_0)), \forall k\in \mathbb{N}.$$ Thus, we obtain that there exists $k'_0\in \mathbb{N}$ such that $$u(\bm{\phi}_{\bm{y}_0}^{\pi'}(k'_0))<0,$$ contradicting the fact that $u(\bm{\phi}_{\bm{y}_0}^{\pi'}(k))\geq 0, \forall k\in \mathbb{N}.$
Therefore, every possible trajectory initialized in the set $\{\bm{x}\in B(\bm{0},R)\mid u(\bm{x})<1\}\setminus X_{\infty}$ will enter the set $X_{\infty}$ in finite time. Consequently, every possible trajectory initialized in the set $\{\bm{x}\in B(\bm{0},R)\mid u(\bm{x})<1\}$ will  approach the equilibrium state $\bm{0}$. 

Combining above arguments, we conclude that $\{\bm{x}\in B(\bm{0},R)\mid u(\bm{x})<1\}$ is a robust region of attraction. \qed
\end{pf}

\begin{remark}
\label{asym}
Note that Theorem \ref{inner} still holds if the equilibrium $\bm{0}$ is asymptotically stable rather than uniformly locally exponentially stable. The proof of Theorem \ref{inner} did not require Assumption \ref{app}.
\end{remark}

\subsection{Theoretical Analysis}
\label{CA}
This section shows that there exists a sequence of solutions to the semi-definite program \eqref{sos} such that their strict one sub-level sets inner-approximate the interior of the maximal robust region of attraction in measure under appropriate assumptions.  

\begin{assumption}
\label{ass2}
 One of the polynomials defining the set $D$ is equal to $h_i^D:=\|\bm{d}\|^2-R_D$ for some constant $R_D\geq 0$, $i\in \{1,\ldots,m_d\}$.
\end{assumption}
Assumption \ref{ass2} is without loss of generality since $D$ is compact, and thus redundant constraint of the form $R_D-\|\bm{d}\|^2\geq 0$ can always be added to the description of $D$ for sufficiently large $R_D$.

\begin{lemma}
\label{uniform}
Under Assumptions \ref{app}, \ref{assump2} and \ref{ass2}, there  exists a sequence $(u_{k_i}(\bm{x}))_{i=0}^{\infty}$ such that $u_{k_i}(\bm{x})$ converges from above to $v$ uniformly over $B(\bm{0},R)$, where $u_{k_i}(\bm{x})\in \mathbb{R}_{k_i}[\bm{x}]$ denotes the $u-$component of a solution to the semi-definite program \eqref{sos} and $v$ is the continuous and bounded solution to the Bellman equation \eqref{e-v}. 
\end{lemma}
\begin{pf}
Let \begin{equation}
\label{reach}
\begin{split}
\Omega(B(\bm{0},R))=\{\bm{y}\in \mathbb{R}^n\mid &\bm{y}=\bm{\phi}_{\bm{x}_0}^{\pi}(i), i\in [0,1],\\
&\bm{x}_0\in B(\bm{0},R), \pi\in \mathcal{D}\}.
\end{split}
\end{equation}
 Since $\bm{f}\in \mathbb{R}[\bm{x},\bm{d}]$, and $D$ and $B(\bm{0},R)$ are compact, $\Omega(B(\bm{0},R))$ is bounded and consequently $\overline{\Omega(B(\bm{0},R))}$ is compact.  Moreover, $B(\bm{0},R)\subset \Omega(B(\bm{0},R))$. Let $B=\overline{B(\bm{0},R)\setminus X_{\infty}}$. 
We infer that for every $\epsilon>0$, there exists a continuous function $v_{\epsilon}$ satisfying \eqref{constraint0} and $|v_{\epsilon}-v|\leq \epsilon$. Obviously, $v_{\epsilon}=v+{\epsilon}$ satisfies such requirement since 
\begin{equation}
\label{constraint1}
\begin{split}
&v_{\epsilon}-v_{\epsilon}(\bm{f})-g\cdot(1-v_{\epsilon})\geq \epsilon c_0, \forall \bm{x}\in B,\forall \bm{d}\in {D},\\
&v_{\epsilon}-1\geq \epsilon, \forall \bm{x}\in B(\bm{0},R)\setminus X,\\
&v_{\epsilon}-h^X_j\geq \epsilon, \forall \bm{x}\in \overline{X}, j=1,\ldots,n_{X},
\end{split}
\end{equation}
where $c_0=\inf\{g(\bm{x})\mid \bm{x}\in B\}$. 
 Since $\overline{\Omega(B(\bm{0},R))}$ is compact, according to Stone-Weierstrass theorem \citep{cotter1990}, there exists a polynomial $u_{k_i}$ of sufficiently high degree $k_i$ such that 
 \[0<u_{k_i}-v_{\epsilon}<\frac{\epsilon}{2} c_0, \forall \bm{x}\in \overline{\Omega(B(\bm{0},R))}.\]
  Thus, we have 
  \begin{equation}
  \label{sosos}
\epsilon<u_{k_i}-v<\epsilon+\frac{\epsilon}{2} c_0, \forall \bm{x}\in \overline{\Omega(B(\bm{0},R))}.  
  \end{equation}
   According to the definition of $\Omega(B(\bm{0},R))$, i.e., \eqref{reach}, we have that $\bm{f}(\bm{x},\bm{d}) \in \Omega(B(\bm{0},R))$ holds for $\bm{x}\in B(\bm{0},R)$ and $\bm{d}\in D$. Therefore, 
\[\epsilon<u_{k_i}(\bm{f}(\bm{x},\bm{d}))-v(\bm{f}(\bm{x},\bm{d}))<\epsilon+\frac{\epsilon}{2} c_0\] holds for $\bm{x}\in B(\bm{0},R)$ and $\bm{d}\in D$. It is easy to check that $u_{k_i}$ satisfies
 \begin{equation}
\label{constraint11}
\begin{split}
&u_{k_i}-u_{k_i}(\bm{f})-g\cdot(1-u_{k_i})>0, \forall \bm{x}\in B,\forall \bm{d}\in D,\\
&u_{k_i}-1>0, \forall \bm{x}\in B(\bm{0},R)\setminus X,\\
&u_{k_i}-h_j^X>0, \forall \bm{x}\in \overline{X}, j=1,\ldots,n_{X}.
\end{split}
\end{equation}
From Putinar's Positivstellensatz \citep{putinar93} and arbitrariness of $\epsilon$, we obtain $u_{k_i}(\bm{x})$ converges from above to $v$ uniformly over $B(\bm{0},R)$ with $i$ approaching infinity. \qed
\end{pf}

Finally, we conclude that $\{\bm{x}\in B(\bm{0},R)\mid u_{k_i}(\bm{x})<1\}$ converges to the interior of the maximal robust region of attraction with $i$ approaching infinity. 
\begin{theorem}
Let $u_{k_i}(\bm{x})$ satisfy the condition in Lemma \ref{uniform}. Then the set $\mathcal{R}_{k_i}=\{\bm{x}\in B(\bm{0},R)\mid u_{k_i}(\bm{x})<1\}$ satisfies 
\[\mathcal{R}_{k_i}\subset \mathcal{R}^{\circ}\text{~and}\]
\[\lim_{i\rightarrow \infty}\mu(\mathcal{R}^{\circ}\setminus \mathcal{R}_{k_i})=0.\]
\end{theorem}
\begin{pf}
$\mathcal{R}_{k_i}\subset \mathcal{R}^{\circ}$ is an immediate consequence of Lemma \ref{uniform} since $u_{k_i}\geq v$ over $B(\bm{0},R)$ according to \eqref{sosos}. 

According to Theorem 1 as well as Theorem 3 in \citep{lasserre2015} and Lemma \ref{uniform}, we have $\lim_{i\rightarrow \infty} \mu(\mathcal{R}^{\circ}\setminus \mathcal{R}_{k_i})=0.$\qed
\end{pf}
\section{Illustrative Examples}
\label{ex}
In this section we evaluate the semi-definite programming based method on two examples. The computations were performed on an i7- 7500U 2.70GHz CPU with 32GB RAM running Windows 10.  YALMIP \citep{lofberg2004} and Mosek \citep{mosek2015mosek} were used to implement  \eqref{sos}.

\begin{table}
\begin{center}
\begin{tabular}{|c|c|c|c|c|c|c|c|}
  \hline
  \multirow{1}{*}{}&\multicolumn{6}{|c|}{\texttt{SDP}} \\\hline
    Ex.&$k$&$d_{s_{i}}$&$d_{s_{3,i_1}}$&$d_{s_{i_2,j}}$&$d_{s_{9,l,j}}$&$T$\\\hline
   \multirow{2}{*}{1}&6&8&8&$8$&8&2.10 \\
                      &10&12&12&$12$&12&17.50\\\hline
    \multirow{3}{*}{2}&4&4&4&$4$&4&5.45 \\
                      &6&6&6&$6$&6&24.67\\
                      &8&8&8&$8$&8&316.12
                      \\\hline
   \end{tabular}
   \oomit{\begin{tabular}{|c|c|c|c|c|c|c|c||c|c|c|c|c|}
  \hline
   \multicolumn{4}{|c|}{\texttt{VI}} \\\hline
    Ex.&$\epsilon$& N& M&$T_{\mathtt{VI}}$\\\hline
   \multirow{1}{*}{1}&$10^{-16}$&$10^4$&50& 748.92 \\\hline
   \end{tabular}}
\end{center}
{Table 1. Parameters of our implementations on Example \ref{ex1}.  $k, d_{s_{i}},d_{s_{3,i_1}}, d_{s_{i_2,j}}, d_{s_{9,l,j}}$: degree of polynomials $u, s_{i},s_{3,i_1}, s_{i_2,j}, s_{9,l,j}$ in \eqref{sos}, respectively, $i_1=1,\ldots,m_d$, $i=0,\ldots,2$, $i_2=4,\ldots,8$, $j=1,\ldots,n_{X}$, $l=1,\ldots,n_{X}$; $T$: computation times (seconds).}
\end{table}

\begin{example}
\label{ex1}
Consider the discrete-generation predator-prey model from \citep{halanay2000stability}, 
\begin{equation*}
\begin{cases}
&x(j+1)=0.5x(j)-x(j)y(j),\\
&y(j+1)=-0.5y(j)+(d(j)+1)x(j)y(j),
\end{cases}
\end{equation*}
where $j\in \mathbb{N}$.

In this example  we consider $D=\{d\in \mathbb{R}\mid d^2-0.01\leq 0\}$ and
$X=\{(x,y)\mid  x^2+y^2<1\}$.  The origin $\bm{0}$ for this example is uniformly locally exponentially stable. $g(x,y)=x^2+y^2$, $R=1.6$, $h_0(x,y)=x^2+y^2$ and $h_{\infty}(x,y)=100(x^2+y^2)$ are used to perform computations on the semi-definite program \eqref{sos}.

The function $h_{\infty}(\bm{x})=100x^2+100y^2$ defining $X_{\infty}$ is a Lyapunov function such that $X_{\infty}\subset X$ is a robust region of attraction. This argument can be justified by first encoding the following constraint $$h_{\infty}(\bm{x})-h_{\infty}(\bm{f}(\bm{x},\bm{d}))>0, \forall \bm{x}\in X_{\infty}\setminus \{\bm{0}\}, \forall \bm{d}\in D$$ in the form of sum-of-squares constraints and then verifying the feasibility of the constructed sum-of-squares constraints, where $\bm{f}(\bm{x},\bm{d})=(0.5x-xy;-0.5y+(d+1)xy)$. Assumption 2(a) is satisfied. Also, the set $B(\bm{0},R)=\{\bm{x}\mid h_0(x,y)\leq 1.6\}$ satisfies Assumption \ref{assump2}(b). Since $X\subseteq B(\bm{0},R)$,  we just need to verify $\{\bm{x}\mid \bm{x}=\bm{f}(\bm{x}_0,\bm{d}), \bm{x}_0\in X, \bm{d}\in D\}$. This argument is justified by first encoding the following constraint \[
\begin{split}
1.6-\big(&0.5x-x y\big)^2\\
&-\big(-0.5y+(d+1)xy\big)^2\geq 0, \forall (x,y)\in X, \forall d\in D
\end{split}\] in the form of sum-of-squares constraints and then verifying its feasibility. Moreover, the function $d^2-0.01$ defining $D$ satisfies Assumption \ref{ass2}.  Therefore,  Lemma \ref{uniform} holds, implying that the existence of solutions to the semi-definite program \eqref{sos} is guaranteed. 

Robust regions of attraction, which are computed via solving the semi-definite program \eqref{sos} with approximating polynomials of degree $6$ and $10$ respectively, are illustrated in Fig. \ref{fig-one-3}.  We observe from Fig. \ref{fig-one-3} that the robust region of attraction computed when $k=10$ approximates the maximal robust region of attraction tightly by comparing with the maximal one estimated via simulation methods.
Here the simulation method requires gridding the state space and the disturbance space, and the check of whether grid states will hit the region $X_{\infty}$ while remaining inside the set $X$ preceding the hitting time. Two trajectories, one respecting the state constraint and one violating the state constraint, are illustrated in Fig. \ref{fig-one-31}. They are generated by extracting the perturbation input $d(j)$ from $D$ randomly for $j\in \mathbb{N}$.

\begin{figure}
\center
\includegraphics[width=3.0in,height=1.3in]{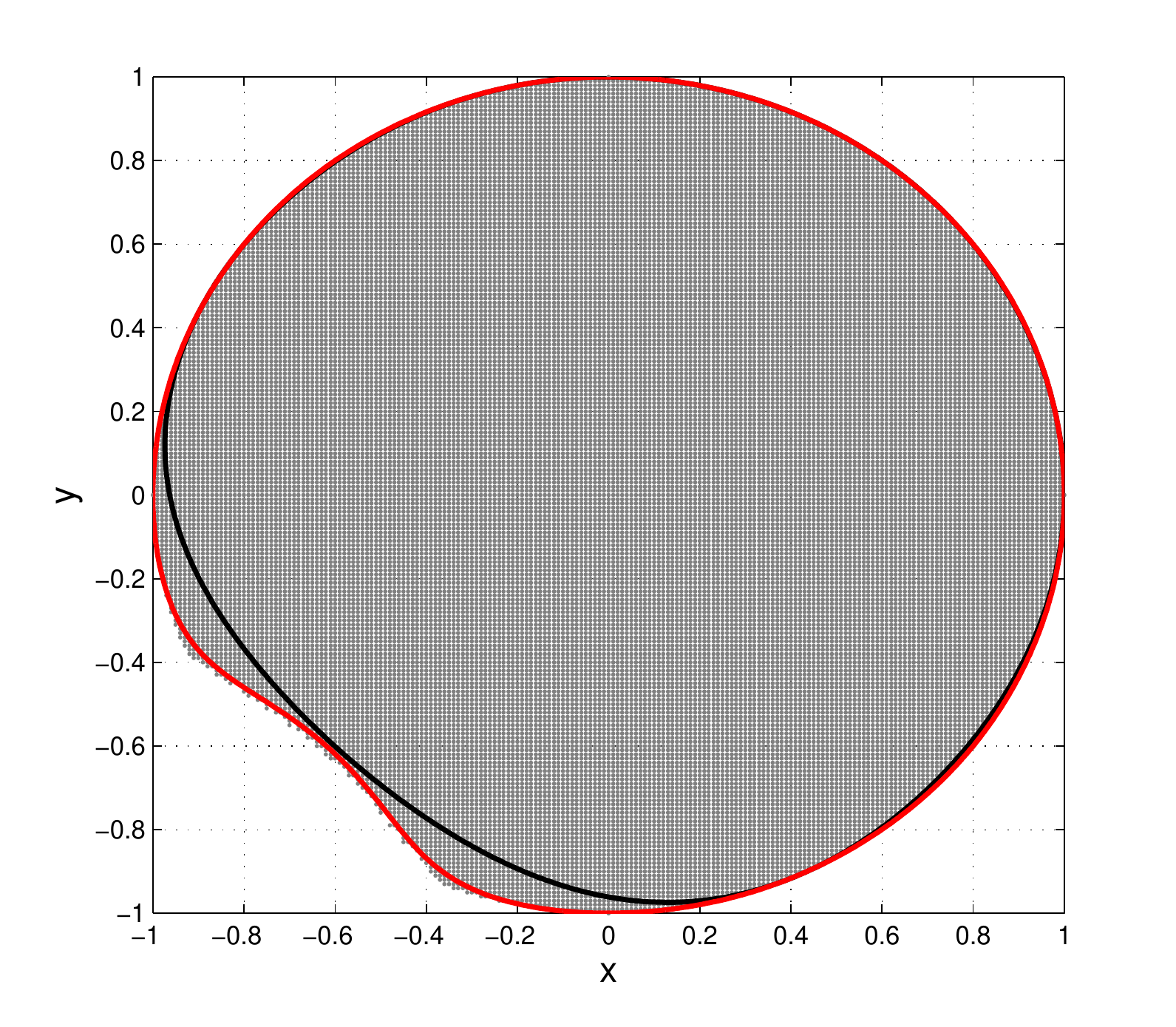} 
\caption{Estimations of $\mathcal{R}$ for Example \ref{ex1}. Black and red curves  denote the boundaries of robust regions of attraction computed when $k=6$ and $k=10$, respectively. Gray region denotes the maximal robust region of attraction estimated via simulation techniques.}
\label{fig-one-3}
\end{figure}

\begin{figure}
\center
\includegraphics[width=3.0in,height=1.3in]{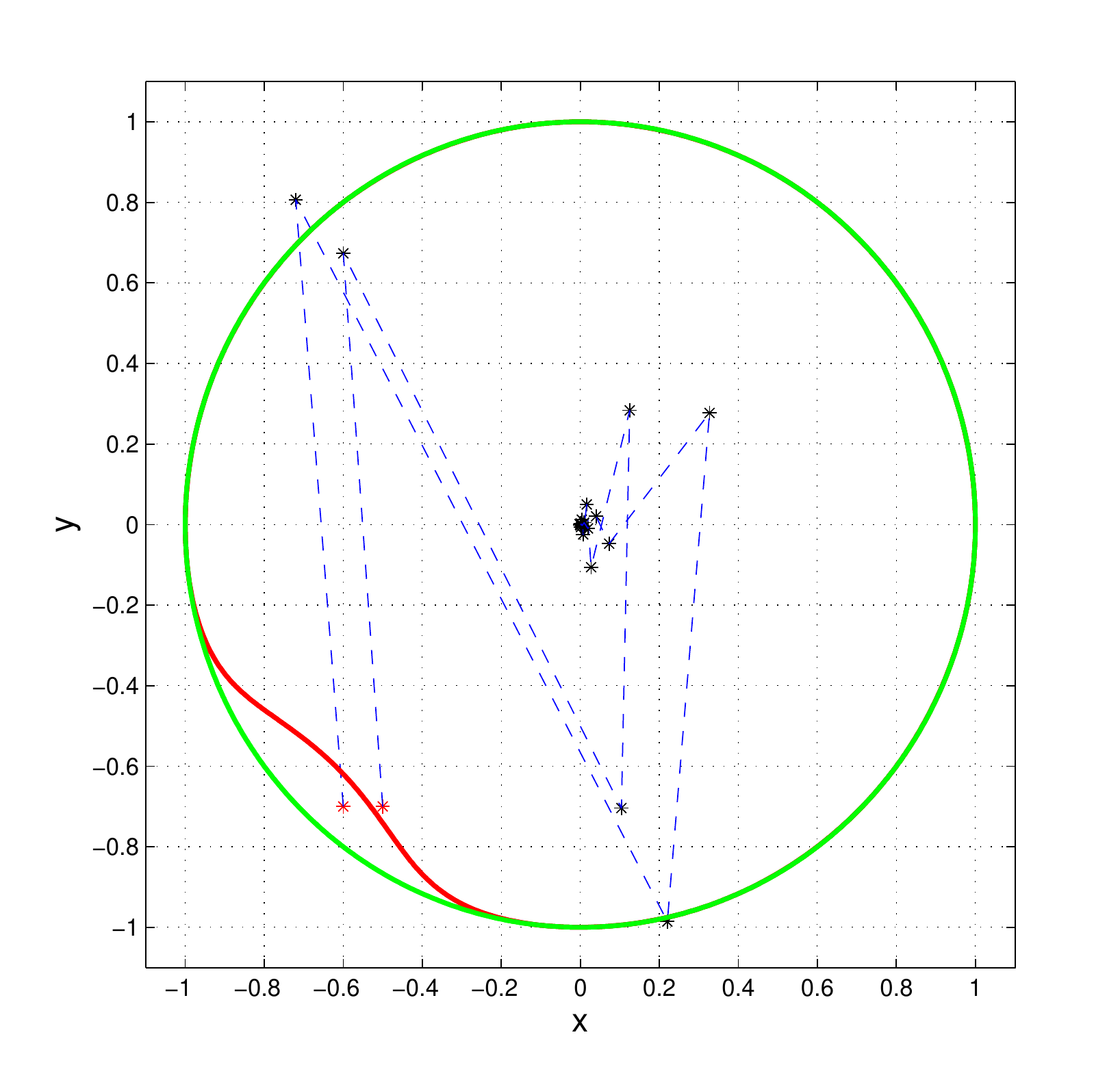} 
\caption{An illustration of two trajectories for Example \ref{ex1}.  Green and red curves  denote the boundaries of the state constraint set $X$ and the robust region of attraction computed when $k=10$. Red stars and black stars denote the initial states  and subsequent states, respectively. The dash blue line denotes the transition between states.}
\label{fig-one-31}
\end{figure}
\end{example}
\begin{example}
\label{ex_three}
Consider a three-dimensional perturbed discrete-time Lotka-Volterra model adopted from \citep{bischi2010three}, 
\begin{equation*}
\begin{cases}
&x(j+1)=x(j)(e_1+a_{1}x(j)+a_{2}y(j)+a_{3}z(j)),\\
&y(j+1)=y(j)(e_2+a_{4}x(j)+a_{5}y(j)+a_{6}z(j)),\\
&z(j+1)=z(j)(e_3+a_{7}x(j)+a_{8}y(j)+a_{9}z(j)),
\end{cases}
\end{equation*}
where $e_1=e_2=e_3=0.5$, $a_{1}=0.5+d$, $a_2=a_6=-0.5$, $a_3=a_4=a_5=a_7=a_8=a_9=0.5$, $D=\{d\in \mathbb{R}\mid d^2-0.01\leq 0\}$ and $X=\{(x,y,z)\in \mathbb{R}^3\mid x^2+y^2+z^2<1\}$. The origin $\bm{0}$ is uniformly locally exponentially stable.

For this example, the sets $X_{\infty}=\{\bm{x}\in \mathbb{R}^3\mid 100(x^2+y^2+z^2)<1\}$ and $B(\bm{0},R)=\{\bm{x}\in \mathbb{R}^3\mid x^2+y^2+z^2\leq 1.6\}$ satisfy Assumption \ref{assump2} and are used for perform computations. Moreover, the function $d^2-0.01$ defining $D$ satisfies Assumption \ref{ass2}. Therefore, Lemma \ref{uniform} holds, implying that the existence of solutions to the semi-definite program \eqref{sos} with $g(\bm{x})=x^2+y^2+z^2$ is guaranteed.

 Plots of computed robust regions of attraction for approximating polynomials of degree $k = 4, 6, 8$ on planes $y-z$ with $x=0$, $x-z$ with $y=0$ and $x-y$ with $z=0$ are shown in Fig. \ref{fig-two-3}. In order to shed light on the accuracy of the computed regions of attraction, we also use
the simulation technique to synthesize estimations of the maximal robust region of attraction on planes $y-z$ with $x=0$, $x-z$ with $y=0$ and $x-y$ with $z=0$ by taking initial states in the state
spaces $\{\bm{x}\in \mathbb{R}^3\mid x^2+y^2+z^2\leq 1\wedge x=0\} $,
$\{\bm{x}\in \mathbb{R}^3\mid x^2+y^2+z^2\leq 1\wedge y=0\}$ and 
$\{\bm{x}\in \mathbb{R}^3\mid x^2+y^2+z^2\leq 1\wedge z=0\}$, respectively.
They are the gray regions in Fig. 5. We observe from Fig. \ref{fig-two-3} that the robust region of attraction computed when $k=8$ could approximate the maximal robust region of attraction well.

\begin{figure}
\center
\includegraphics[width=3.0in,height=1.3in]{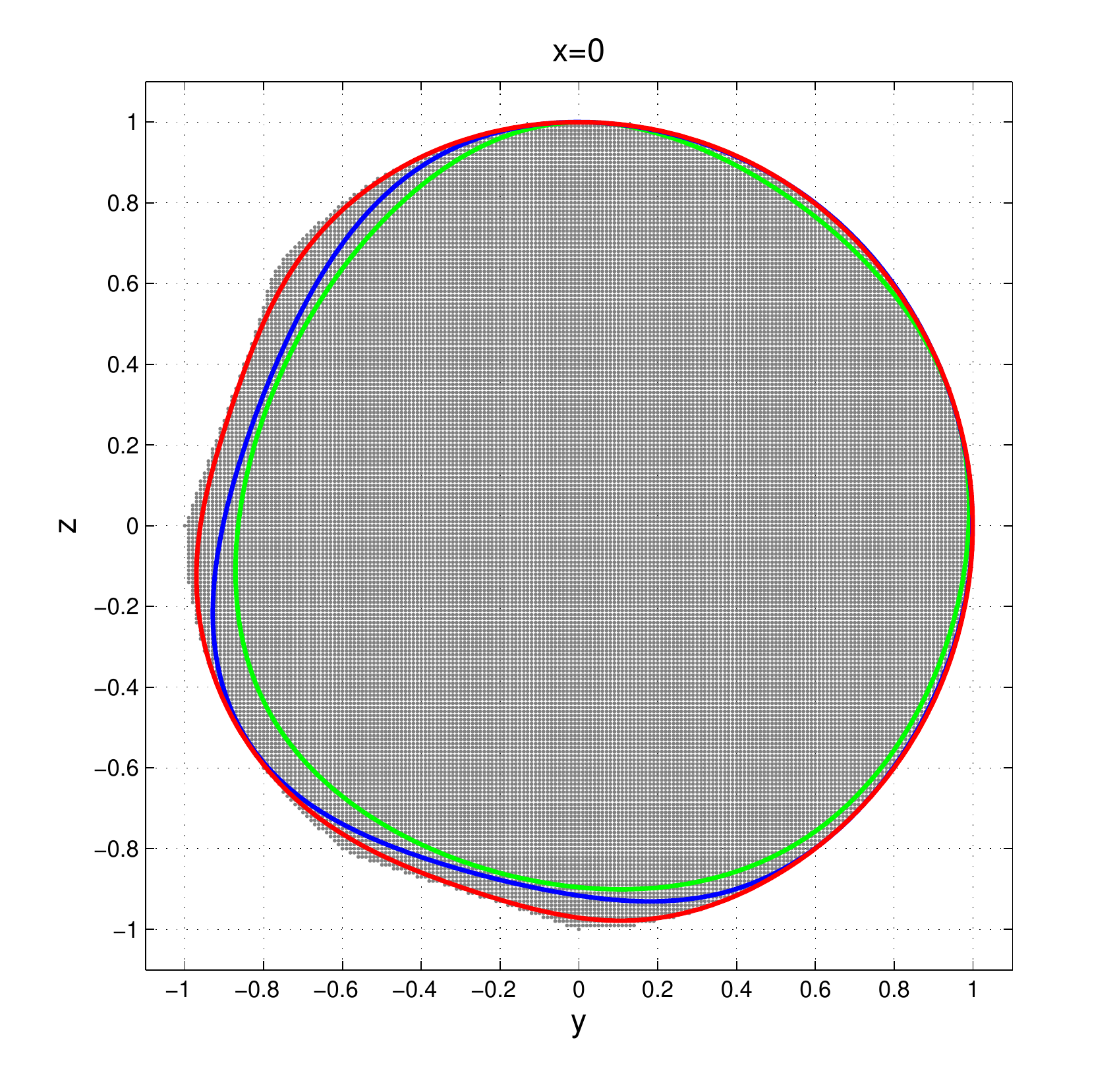} \\
\includegraphics[width=3.0in,height=1.3in]{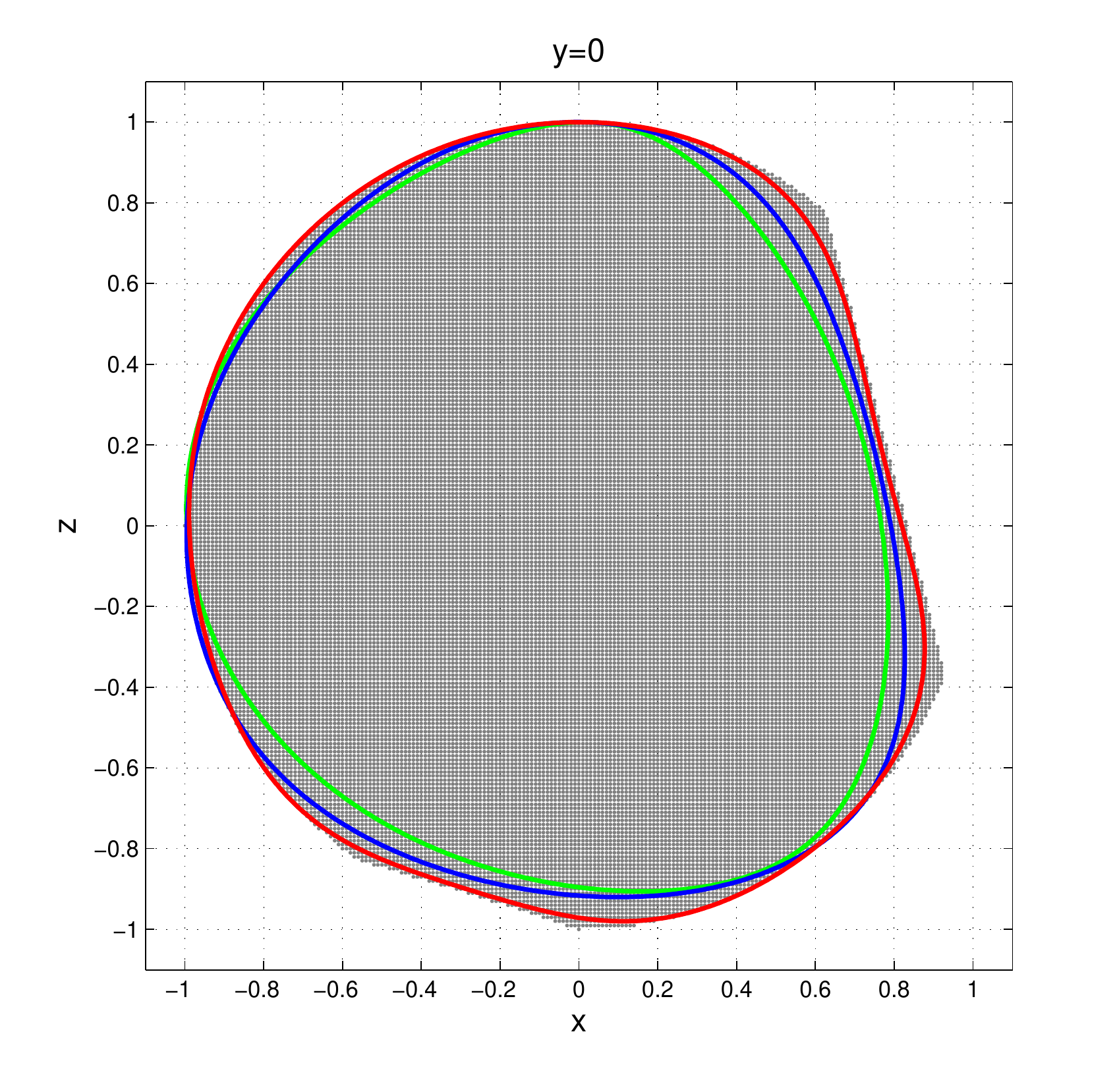} \\
\includegraphics[width=3.0in,height=1.3in]{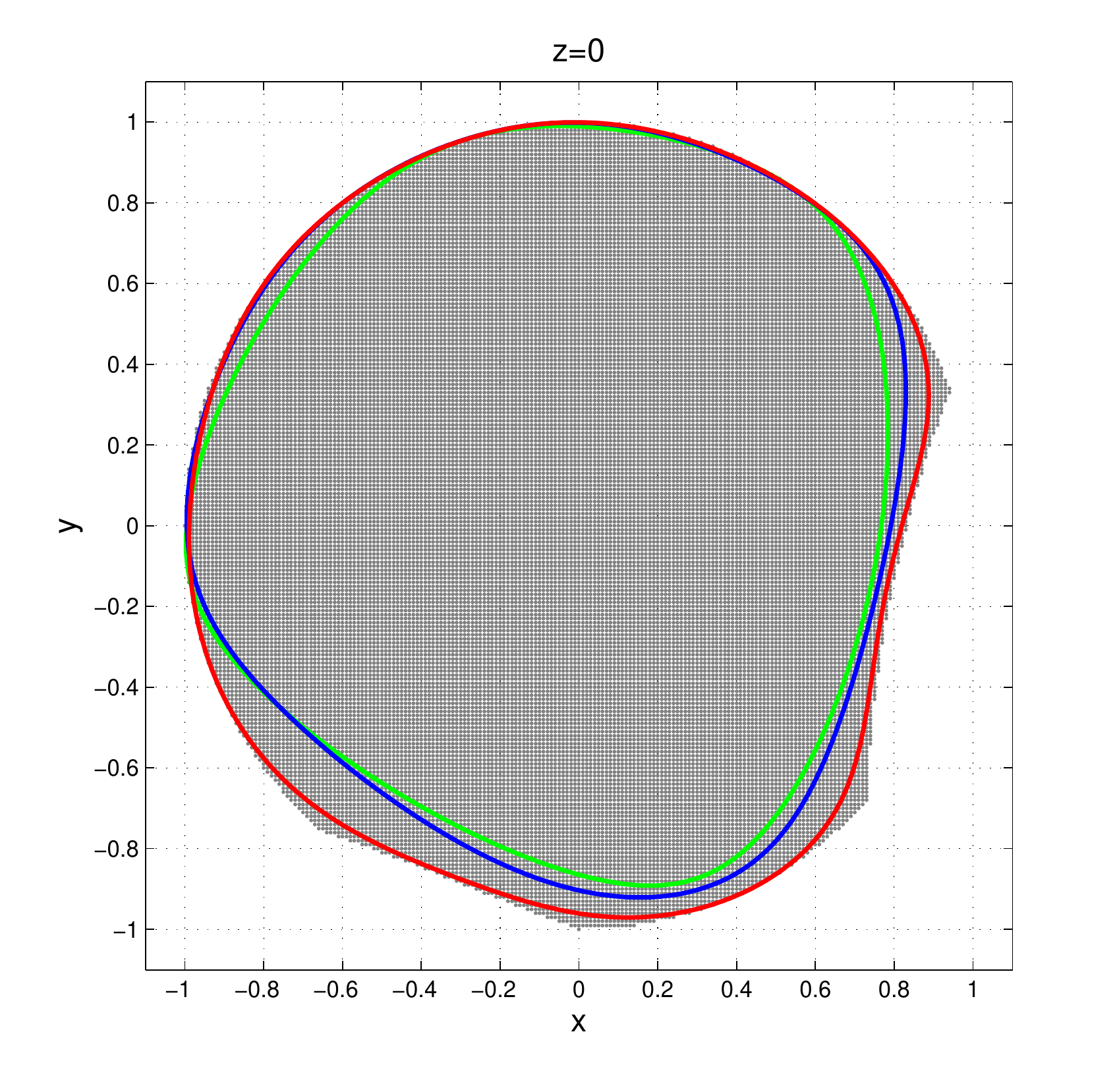} 
\caption{Estimations of $\mathcal{R}$ for Example \ref{ex_three}. Green, blue and red curves denote the boundaries of robust regions of attraction computed when $k=4$, $6$ and $8$, respectively. Gray region denotes the maximal robust region of attraction estimated via simulation techniques.}
\label{fig-two-3}
\end{figure}

\end{example}

\section{Conclusion}
\label{con}
In this paper we proposed  a semi-definite programming based  method for computing robust regions of attraction for state-constrained perturbed discrete-time polynomial  systems. The semi-definite program was constructed based on a Bellman equation. Also, there exists a sequence of solutions to the semi-definite program such that their strict one sub-level sets inner-approximate the interior of the maximal robust region of attraction in measure under appropriate assumptions. Two examples demonstrated the performance of our approach. 

In near future we would like to compare the proposed method in this paper with existing methods on estimating robust regions of attraction for discrete-time systems. Also, we would extend the proposed method for computing robust regions of attraction of state-constrained perturbed continuous-time polynomial systems.

\bibliography{ifacconf} 
\section{Appendix}
In this section we characterize the interior of the maximal robust region of attraction  $\mathcal{R}$ as the strict one sub-level set of the unique bounded and continuous solution to the modified Bellman equation \eqref{e-v}. Its derivation follows that in \citep{grune2015}.

\subsection{Robust Regions of Uniform Attraction}
\label{UV}
In this subsection we introduce the maximal robust region of uniform attraction, which is equal to the interior of the maximal robust region of attraction. The maximal robust region of uniform attraction was first proposed in \citep{grune2015} for state-constrained perturbed continuous-time systems.  

Denote the first hitting time $k'(\bm{x}_0,\pi)$, induced by the initial state $\bm{x}_0$ and the input policy $\pi$, of $B(\bm{0},\overline{\epsilon})$ as
\begin{equation}
    \label{hitting}
k'(\bm{x}_0,\pi):=\inf\{k>0\mid \bm{\phi}_{\bm{x}_0}^{\pi}(k)\in B(\bm{0},\overline{\epsilon})\},
\end{equation}
where $B(\bm{0},\overline{\epsilon})$ is defined in \eqref{half}.  Also, let the Euclidean distance between a point $\bm{x}\in \mathbb{R}^n$ and a set $A\subset \mathbb{R}^n$ be $\mathtt{dist}(\bm{x},A):=\mathtt{inf}_{\bm{y}\in A}\|\bm{x}-\bm{y}\|$, and the set of $\delta$-admissible perturbation input policies be 
 \[\mathcal{D}_{ad,\delta}(\bm{x}_0):=\{\pi\mid \mathtt{dist}(\bm{\phi}_{\bm{x}_0}^{\pi}(k),X^{c})>\delta \text{ for } k\in \mathbb{N}\},\]
 where $\delta>0$ and $X^c$ is the complement of the set $X$. The maximal robust region of uniform attraction $\mathcal{R}_0$ is then defined by 
\begin{equation*}
\mathcal{R}_0: =\left\{\bm{x}_0\in \mathbb{R}^n\middle|\;
\begin{aligned}
&\text{there exists~} \delta>0 \text{~s.t.~}\mathcal{D}_{ad,\delta}(\bm{x}_0)\\
&=\mathcal{D}\text{~and~}\sup_{\pi\in \mathcal{D}}k'(\bm{x}_0,\pi)<\infty
\end{aligned}
\right\}.
\end{equation*}

\oomit{Before presenting the openness property of the region $\mathcal{R}_0$ and the relationship between $\mathcal{R}_0$ and $\mathcal{R}$, 
we first characterize the maximal robust region of attraction $\mathcal{R}$ using $k'(\bm{x}_0,\pi)$.
\begin{lemma}
\label{1}
Under Assumption \ref{app}, the maximal robust region of attraction $\mathcal{R}$ is equal to $\mathcal{R}'$, where 
\begin{equation*}
\mathcal{R}'=\left\{\bm{x}_0 \middle|\;
\begin{aligned}
&\mathcal{D}_{ad}(\bm{x}_0)=\mathcal{D}\text{ and }\\
& k'(\bm{x}_0,\pi)<\infty \text{ for }\pi\in \mathcal{D}
\end{aligned}
\right\}.
\end{equation*}
\end{lemma} 
\begin{pf}
Obviously, $\mathcal{R}'\subseteq \mathcal{R}$. Thus, we just prove $\mathcal{R}\subseteq \mathcal{R}'$.

Assume that there exists $\bm{y}_0\in \mathcal{R}\setminus \mathcal{R}'$. Obviously, $\bm{\phi}_{\bm{y}_0}^{\pi}(l)\in X$ for $l\in \mathbb{N}$ and $\pi \in \mathcal{D}$. Moreover, there exists a perturbation input policy $\pi'$ such that $\bm{\phi}_{\bm{y}_0}^{\pi'}(l)\notin B(\bm{0},\overline{\epsilon})$ for $l\in \mathbb{N}$, contradicting $\lim_{l\rightarrow \infty}\bm{\phi}_{\bm{y}_0}^{\pi}(l)=\bm{0}$ for $\pi\in \mathcal{D}$. Therefore, $\mathcal{R}\subseteq \mathcal{R}'$. 

Above all, we have $\mathcal{R}=\mathcal{R}'$. \qed
\end{pf}}

 Lemma \ref{lemma1} presents the openness property of the region $\mathcal{R}_0$ as well as the relationship between $\mathcal{R}_0$ and $\mathcal{R}$.
\begin{lemma}
\label{lemma1}
Under Assumption \ref{app}, then
\begin{enumerate}
\item[(a)] $\mathcal{R}_0=\mathcal{R}'_0$, where 
\begin{equation*}
\label{R2}
\mathcal{R}'_0=\left\{\bm{x}_0 \in \mathbb{R}^n \middle|\;
\begin{aligned}
&\text{there exists~} \delta>0 \text{ s.t. }\mathcal{D}_{ad,\delta}(\bm{x}_0)\\
&=\mathcal{D}
\text{~and there exists a function~} \\
&\beta(k): \mathbb{N}\rightarrow [0,\infty)\text{~satisfying~}\\&\lim_{k\rightarrow \infty} \beta(k)=0 \text{~and~}\\
&\|\bm{\phi}_{\bm{x}_0}^{\pi}(k)\|\leq \beta(k)\text{~for~} k\in \mathbb{N}\\
& \text{~and~} \pi\in \mathcal{D}.
\end{aligned}
\right\}.
\end{equation*}
\item[(b)] $\mathcal{R}_0$ is open.
\item[(c)] $\mathcal{R}_0=\mathcal{R}^{\circ}$.
\end{enumerate}
\end{lemma}
\begin{pf}
(a). Let $\bm{x}_0\in \mathcal{R}_0$ and $K=\sup_{\pi\in \mathcal{D}}k'(\bm{x}_0,\pi)<\infty$. Then, for $k\geq K$ we have $$\|\bm{\phi}_{\bm{x}_0}^{\pi}(k)\|\leq \beta(r,k)=\lambda^kM r,$$ where $r$ is defined in \eqref{r}. Hence, for $k\geq K$ we can choose $\beta(k)=\beta(r,k).$ Since $\bm{\phi}_{\bm{x}_0}^{\pi}(k) \in X$ for $k\in [0,K]$ and $\pi \in \mathcal{D}$, and $X$ is bounded, there exists $M'\geq 0$ such that $$\|\bm{\phi}_{\bm{x}_0}^{\pi}(k)\|\leq M', \forall k\in [0,K], \forall \pi \in \mathcal{D}.$$ Choosing $\beta(k)=M'$ for $k\in [0,K]$ then yields the function $\beta(k)$ with the desired properties. Thus, $$\bm{x}_0\in \mathcal{R}'_0,$$ implying that $\mathcal{R}_0\subseteq \mathcal{R}'_0.$

Conversely, let $\bm{x}_0\in \mathcal{R}'_0$ and pick the corresponding $\delta>0$ and $\beta(k)$. Then there exists $K>0$ such that
 $$\beta(k)<\overline{\epsilon}, \forall k\geq K$$ ($K$ exists since $\lim_{k\rightarrow \infty}\beta(k)=0$), where $\overline{\epsilon}$ is defined in \eqref{half}. Then we have $$\|\bm{\phi}_{\bm{x}_0}^{\pi}(k)\|\leq \beta(k)<\overline{\epsilon}, \forall k\geq K, \forall \pi \in \mathcal{D},$$ which implies $$\bm{\phi}_{\bm{x}_0}^{\pi}(k)\in B(\bm{0},\overline{\epsilon}), \forall k\geq K, \forall \pi\in \mathcal{D}.$$ Hence, $$k'(\bm{x}_0,\pi)\leq K, \forall \pi\in \mathcal{D}$$ and thus $$\sup_{\pi\in \mathcal{D}}k'(\bm{x}_0,\pi)\leq K<\infty.$$ Also, since $\mathcal{D}_{ad,\delta}=\mathcal{D}$, we have that $\bm{x}_0\in \mathcal{R}_0,$ implying that $\mathcal{R}'_0\subseteq \mathcal{R}_0.$

(b). Since $\mathcal{R}_0=\mathcal{R}'_0$, we prove the openness of $\mathcal{R}'_0$ instead. Let $\bm{x}_0\in \mathcal{R}'_0$ with  corresponding $\delta>0$ and $\beta(\cdot):\mathbb{N}\rightarrow [0,\infty)$, and $K>0$ be such that $\beta(k)<\frac{\overline{\epsilon}}{2}$ for $k\geq K$, where $\overline{\epsilon}$ is defined in \eqref{hitting}.

Since $\bm{f}(\bm{x},\bm{d})$ is Lipschitz continuous over $\bm{x}\in X$ uniformly over $\bm{d}\in D$, implying that there exists $B(\bm{x}_0,\epsilon)$ such that for $\bm{y}_0\in B(\bm{x}_0,\epsilon), \pi\in \mathcal{D}$ and $k\in [0,K]$, \[\|\bm{\phi}_{\bm{x}_0}^{\pi}(k)-\bm{\phi}_{\bm{y}_0}^{\pi}(k)\|<\min\{\frac{\delta}{2},\frac{\overline{\epsilon}}{2}\}.\] This further implies that for $\bm{y}_0\in B(\bm{x}_0,\epsilon)$, $\pi\in \mathcal{D}$ and $k\in [0,K]$, $$\mathtt{dist}(\bm{\phi}_{\bm{y}_0}^{\pi}(k),X^c)>\frac{\delta}{2}$$ holds. Thus, $\bm{\phi}_{\bm{y}_0}^{\pi}(K)\in B(\bm{0},\overline{\epsilon}), \forall \pi\in \mathcal{D}.$ Hence $$\sup_{ \pi\in \mathcal{D}}k'(\bm{y}_0,\pi)\leq K.$$ Together with \eqref{half} this implies $$\mathcal{D}_{ad,\min\{\frac{\delta}{2},\frac{r}{2}\}}(\bm{y}_0)=\mathcal{D},$$ hence we conclude that $\bm{y}_0\in \mathcal{R}'_0$. Thus, $B(\bm{x}_0,\epsilon)\subset \mathcal{R}'_0$ and consequently $\mathcal{R}'_0$ is open.

(c). Obviously, $\mathcal{R}_0\subseteq \mathcal{R}$. Therefore, $\mathcal{R}_0^{\circ}\subseteq \mathcal{R}^{\circ}$ and by (b) it implies $\mathcal{R}_0\subseteq \mathcal{R}^{\circ}$. 

Next we just prove that $\mathcal{R}^{\circ}\subseteq \mathcal{R}_0$. Let $\bm{x}_0\in \mathcal{R}^{\circ}\setminus \mathcal{R}_0$. Since $\bm{x}_0\notin \mathcal{R}_0$,  either 
\begin{equation}
\label{k'}
\sup_{\pi\in \mathcal{D}}k'(\bm{x}_0,\pi)=\infty
\end{equation} or 
\begin{equation}
\label{D}
\mathcal{D}_{ad,\delta}(\bm{x}_0)\neq \mathcal{D}, \forall \delta>0
\end{equation} 
must hold. If \eqref{k'} holds, then we obtain $\bm{x}_0\in \partial \mathcal{R}$ since in every neighborhood of $\bm{x}_0$ there exist $\bm{x}'_0$ and a perturbation input policy $\pi$ such that $k'(\bm{x}'_0,\pi)=\infty$, contradicting $\bm{x}_0\in \mathcal{R}^{\circ}$.

Hence assume $$K=\sup_{\pi\in \mathcal{D}}k'(\bm{x}_0,\pi)<\infty.$$ Then we have the conclusion that \eqref{D} holds and thus there exists a sequence $(\pi_{i},k_i)_{i\in \mathbb{N}}$ such that $$\lim_{i\rightarrow \infty}\mathtt{dist}(\bm{\phi}_{\bm{x}_0}^{\pi_{i}}(k_i),X^c)=0.$$ Since \eqref{half} and $k'(\bm{x}_0,\pi_i)\leq K$, we have that $$k_i\leq K, \forall i\in \mathbb{N}.$$ Also, since $\bm{x}_0\in \mathcal{R}$, we have that $\bm{\phi}_{\bm{x}_0}^{\pi}(j) \in X$ for $j\in \mathbb{N}$ and $\pi \in \mathcal{D}$. Thus,  $\bm{x}_i=\bm{\phi}_{\bm{x}_0}^{\pi_{i}}(k_i)$ is bounded. The fact that $\bm{f}(\bm{x},\bm{d})$ is locally Lipschitz continuous over $\mathbb{R}^n$ yields that for every $\epsilon>0$, the set $$\{\bm{\phi}_{\bm{y}}^{\pi_{i}}(k_i)\mid \bm{y}\in B(\bm{x}_0,\epsilon)\}$$ contains a ball $B(\bm{x}_i,\rho)$ with $\rho>0$ independent of $i$(since $k_i\leq K$, $\forall i\in \mathbb{N}$). For sufficiently large $i$ this implies $B(\bm{x}_i,\rho)\nsubseteq X$. This means that $$\pi_{i}\notin \mathcal{D}_{ad,0}(\bm{z}_i)$$ for some $\bm{z}_i\in B(\bm{x}_0,\epsilon)$ and consequently $\bm{z}_i\notin \mathcal{R}.$ Since $\epsilon>0$ is arbitrary, this implies $\bm{x}_0\in \partial \mathcal{R},$ again contradicting $\bm{x}_0\in \mathcal{R}^{\circ}$. Hence, $\mathcal{R}^{\circ}\setminus \mathcal{R}_0=\emptyset,$ implying $\mathcal{R}^{\circ}\subset \mathcal{R}_0.$ \qed
\end{pf}

\subsection{Bellman Equations}
\label{BEE}
In this section we mainly present a modified Bellman equation, to which the strict one sub-level set of the unique bounded and continuous solution is equal to the maximal robust region of uniform attraction $\mathcal{R}_0$. For this sake we first introduce a value function, whose strict one sub-level set is equal to the maximal robust region of uniform attraction $\mathcal{R}_0$. Then we reduce this value function to the unique continuous and bounded solution to a modified Bellman equation. 

We first introduce a semi-definite positive polynomial cost $g:\mathbb{R}^n\rightarrow \mathbb{R}$ satisfying that $g(\bm{x})=0$ iff $\bm{x}=\bm{0}$. For the sake of simplicity, we denote $\ln(g(\bm{\phi}_{\bm{x}}^{\pi}(i))+1)$ and $\ln(l(1-h_j^X(\bm{\phi}_{\bm{x}}^{\pi}(i))))$
as $g_{i}(\bm{x},\pi)$ and $h_{j,i}(\bm{x},\pi)$ respectively, i.e., 
\[g_{i}(\bm{x},\pi)=\ln(g(\bm{\phi}_{\bm{x}}^{\pi}(i))+1)\]
and  
\begin{equation}
\label{hij}
h_{j,i}(\bm{x},\pi)=\ln(l(1-h_{j}^X(\bm{\phi}_{\bm{x}}^{\pi}(i)))),
\end{equation}
where 
\begin{equation*}
\begin{split}
    l(x)=
    \left\{
                \begin{array}{lll}
                  x,& \text{if } x\geq 0, \\
                  0, &\text{otherwise.}\\                  
                \end{array}
              \right.
              \end{split}
\end{equation*}
Besides, we define $\ln 0:=-\infty$.

We define the value function $V:\mathbb{R}^n\rightarrow \mathbb{R}^{+}\cup \{\infty\}$ as
 \begin{equation}
 \label{V0}
 \begin{split}
 V(\bm{x}):=\sup_{\pi\in \mathcal{D}}&\sup_{k\in \mathbb{N}}\big\{\sum_{i=1}^k g_{i-1}(\bm{x},\pi)-\min_{j\in \{1,\ldots,n_{X}\}}h_{j,k}(\bm{x},\pi)\big\}
 \end{split}
 \end{equation}
and consider the Kruzhkov transformed optimal value function $v:\mathbb{R}^n\rightarrow [0,1]$ given by
\begin{equation}
\label{v}
 v(\bm{x}):=1-e^{-V(\bm{x})}=\sup_{\pi\in \mathcal{D}}\sup_{k\in \mathbb{N}}\big\{1-e^{\tilde{V}}\big\},
\end{equation}
where 
\begin{equation}
\label{VVV}
\begin{split}
\tilde{V}=-\sum_{i=1}^k g_{i-1}(\bm{x},\pi)+
\min_{j\in \{1,\ldots,n_{X}\}}h_{j,k}(\bm{x},\pi).
\end{split}
\end{equation}

\begin{theorem}
\label{conti}
Under Assumption \ref{app}, then
\begin{enumerate}
\item[(a)] $\mathcal{R}_0=\{\bm{x}\in \mathbb{R}^n\mid V(\bm{x})<\infty\}=\{\bm{x}\in \mathbb{R}^n\mid v(\bm{x})<1\}.$
\item[(b)] $V(\bm{x})$ is continuous over $\mathcal{R}_0$. Also,  $V(\bm{x})=\infty$ for $\bm{x}\notin \mathcal{R}_0$. 
\item[(c)] $v(\bm{x})$ is continuous over $\mathbb{R}^n$.
\end{enumerate} 
\end{theorem}
\begin{pf}
 In these proofs, $\Omega(\bm{x}_0, k)$ denotes the set of states visited by system \eqref{systems} initialized at $\bm{x}_0$ within $k\geq 1$ steps, i.e., $\Omega(\bm{x}_0, k)=\{\bm{y}\in \mathbb{R}^n\mid \bm{y}=\bm{\phi}_{\bm{x}_0}^{\pi}(i), \forall i\in [0,k]\cap \mathbb{N}, \forall \pi\in \mathcal{D}\}$.  

(a). Firstly, by \eqref{v}, we obtain immediately the equality between the two sets $\{\bm{x}\in \mathbb{R}^n\mid V(\bm{x})<\infty\}$ and $\{\bm{x}\in \mathbb{R}^n\mid v(\bm{x})<1\}$. It remains to prove the first identity that $\mathcal{R}_0=\{\bm{x}\in \mathbb{R}^n\mid V(\bm{x})<\infty\}$.

Let $\bm{x}_0\in \mathcal{R}_0$. We first prove that $$\sup_{\pi\in \mathcal{D}} \sum_{i=1}^{\infty} g_{i-1}(\bm{x}_0,\pi)<\infty.$$

Let $W(\bm{x}_0)=\sup_{\pi\in \mathcal{D}} \sum_{i=1}^{\infty} g_{i-1}(\bm{x}_0,\pi)$.  According to Assumption \ref{app} and the definition of $\mathcal{R}_0$, there exists $K>0$ such that $\bm{\phi}_{\bm{x}_0}^{\pi}(k)\in B(\bm{0},r)$ for $k\geq K$ and $\pi\in \mathcal{D}$. Also, the closure of the reachable set $\overline{\Omega(\bm{x}_0,K)}$ is compact. Thus for $\pi\in \mathcal{D}$,
\begin{equation*}
\begin{split}
W(\bm{x}_0)\leq &K\sup_{\pi\in \mathcal{D},\bm{x}\in \overline{\Omega(\bm{x}_0,K)}}\ln(g(\bm{x})+1)\\
&+\sum_{i=K+1}^{\infty}L_rMr\lambda^{i-K-1}\leq C,
\end{split}
\end{equation*}
 where $L_r$ is the Lipschitz constant of $\ln(g(\bm{x})+1)$ over $\bm{x}\in B(\bm{0},r)$. Therefore $W(\bm{x}_0)<\infty$. Next we prove that 
\[-\sup_{\pi\in \mathcal{D},k\in \mathbb{N}}\min_{j\in \{1,\ldots,n_{X}\}}h_{j,k}(\bm{x}_0,\pi)<\infty.\]
Since $\|\bm{\phi}_{\bm{x}_0}^{\pi}(k)\|\leq \beta(k)$ for $\pi\in \mathcal{D}$, the reachable set $\Omega(\bm{x}_0,\infty)$ is bounded, hence $\overline{\Omega(\bm{x}_0,\infty)}$ is compact. Moreover, since $\mathcal{D}=\mathcal{D}_{ad,\delta}(\bm{x}_0)$ for some $\delta>0$, we have that $\overline{\Omega(\bm{x}_0,\infty)}\subset X$. Also, since each $h^X_j$, $j=1,\ldots, n_X$, is continuous over $X$, it will attain a finite maximum being less than 1 on $\overline{\Omega(\bm{x}_0,\infty)}$ and thus $$\sup_{\pi\in D,k\in \mathbb{N}}\min_{j\in \{1,\ldots,n_{X}\}}h_{j,k}(\bm{x}_0,\pi)$$ will attain a finite minimum over $\overline{\Omega(\bm{x}_0,\infty)}$ according to \eqref{hij}. We prove the claim. 

Let $\bm{x}_0\notin \mathcal{R}_0$. Then either $\sup_{\pi\in \mathcal{D}} k'(\bm{x}_0,\pi)=\infty$ or the existence of $\delta$ in the definition of $\mathcal{R}_0$ is not satisfied, where $k'(\bm{x}_0,\pi)$ is defined in \eqref{hitting}. 

For the first case, there exists a sequence $(\pi_{j'}\in \mathcal{D})_{j'\in \mathbb{N}}$ such that $\lim_{j'\rightarrow \infty}k'(\bm{x}_0,\pi_{{j'}})=\infty$. Then for any $j'\in \mathbb{N}$,
\begin{equation*}
\begin{split}
\sum_{i=1}^{\infty}g_{i-1}(\bm{x}_0,\pi_{{j'}})&\geq \sum_{i=1}^{k'(\bm{x}_0,\pi_{{j'}})}g_{i-1}(\bm{x}_0,\pi_{{j'}})\\
&\geq \ln(c_0+1)k'(\bm{x}_0,\pi_{{j'}}),
\end{split}
\end{equation*}
where $c_0$ is a constant such that $\inf_{\bm{x}\notin B(\bm{0},r)}g(\bm{x})\geq c_0$ (Such $c_0$ exists since $g(\bm{x})$ is a polynomial function over $\bm{x}$ and $g(\bm{x})> 0$ for $\bm{x}\neq \bm{0}$). It follows that $W(\bm{x}_0)\geq \lim_{j'\rightarrow \infty}\sum_{i=1}^{\infty}g_{i-1}(\bm{x}_0,\pi_{{j'}})=\infty$. Therefore, 
$V(\bm{x}_0)=\infty$ since $V(\bm{x}_0)\geq W(\bm{x}_0)$.  In the second case, the non-existence of $\delta$ implies the existence of a sequence $(\pi_{{j'}},k_{j'})_{j'\in \mathbb{N}}$ with $\lim_{j'\rightarrow \infty}\mathtt{dist}(\bm{\phi}_{\bm{x}_0}^{\pi_{{j'}}}(k_{j'}),X^c)= 0$. Then either there exists $l_0\in \mathbb{N}$ such that $\bm{\phi}_{\bm{x}_0}^{\pi_{{l_0}}}(k_{l_0})\in X^c$ or there exists a subsequence $(\bm{x}_{k_{j'_l}})_{l\in \mathbb{N}}$ converging to some $\bm{x}\notin X$ (This is due to the fact that the sequence $(\bm{\phi}_{\bm{x}_0}^{\pi_{{j'}}}(k_{j'}))_{j'\in \mathbb{N}}$ lies in the bounded set $X$.), where $\bm{x}_{k_{j'_l}}=\bm{\phi}_{\bm{x}_0}^{\pi_{{j'_l}}}(k_{j'_l})$. Both cases imply that $$\lim_{j'\rightarrow \infty}\sup_{\pi\in \mathcal{D}}\big(-\min_{j\in \{1,\ldots,n_{X}\}} h_{j,k_{j'}}(\bm{x}_{0},\pi)\big)= \infty.$$ Also, since $$V(\bm{x}_0)\geq \sup_{\pi\in \mathcal{D}}\sup_{j'\in \mathbb{N}}(-\min_{j\in \{1,\ldots,n_{X}\}} h_{j,k_{j'}}(\bm{x}_0,\pi)),$$ we obtain 
$V(\bm{x}_0)=\infty$.

(b). Let $\bm{x}_0, \bm{y}_0 \in \mathcal{R}_0$, 
\begin{equation*}
\begin{split}
|V(\bm{x}_0)-V(\bm{y}_0)|\leq |W(\bm{x}_0)-W(\bm{y}_0)|+|W'(\bm{x}_0)-W'(\bm{y}_0)|,
\end{split}
\end{equation*}
where $W(\bm{x}_0)=\sup_{\pi\in \mathcal{D}} \sum_{i=1}^{\infty} g_{i-1}(\bm{x}_0,\pi)$ and $W'(\bm{x}_0)=\sup_{\pi\in \mathcal{D}}\sup_{k\in \mathbb{N}}\min_{j\in \{1,\ldots,n_{X}\}} h_{j,k}(\bm{x}_0,\pi)$. In the following we separately prove the continuity of $W(\bm{x}_0)$ and $W'(\bm{x}_0).$  Firstly, we prove that $W$ is continuous on $B(\bm{0},\frac{r}{M})$. Assume that $\bm{x}_0\in B(\bm{0},\frac{r}{M})$. Then 
\begin{equation*}
\begin{split}
\sum_{i=0}^{\infty}|\ln(g(\bm{\phi}_{\bm{x}_0}^{\pi}(i))&+1)|\leq L_r \sum_{i=0}^{\infty}\|\bm{\phi}_{\bm{x}_0}^{\pi}(i)\|\\
&\leq L_r M\sum_{i=0}^{\infty} \lambda^i\|\bm{x}_0\|\leq M_1 \|\bm{x}_0\|,
\end{split}
\end{equation*}
where $L_r$ is the Lipschitz constant of $\ln(g(\bm{x})+1)$ over $\bm{x}\in B(\bm{0},r)$, $r$, $\lambda$ and $M$ are defined in \eqref{r}.

For arbitrary but fixed $\epsilon>0$, we can conclude from Assumption \ref{app} that there exists $K>0$ such that 
$M_1\|\bm{\phi}_{\bm{x}_0}^{\pi}(k)\|\leq \frac{\epsilon}{3}$ for $k\geq K$ and $\bm{x}_0\in B(\bm{0},\frac{r}{M})$. In addition, by Lipschitz continuity of $\bm{f}$ there exists $\delta'>0$ such that \[\|\bm{\phi}_{\bm{x}_0}^{\pi}(k)-\bm{\phi}_{\bm{y}_0}^{\pi}(k)\|\leq \frac{\epsilon}{3L_r(K+1)}\] for $k\in [0,K]$ and $\bm{y}_0\in \{\bm{x}\in B(\bm{0},\frac{r}{M})\mid \|\bm{x}-\bm{x}_0\|<\delta'\}$. Then, we have 
\begin{equation}
\label{rr1}
\begin{split}
&|W(\bm{x}_0)-W(\bm{y}_0)|\\
&\leq \sup_{\pi\in \mathcal{D}}\sum_{i=1}^{\infty}|\ln(g(\bm{\phi}_{\bm{x}_0}^{\pi}(i-1))+1)-\ln(g(\bm{\phi}_{\bm{y}_0}^{\pi}(i-1))+1)|\\
&\leq \sup_{\pi\in \mathcal{D}}\big(\sum_{i=0}^{K}L_r\|\bm{\phi}_{\bm{x}_0}^{\pi}(i)-\bm{\phi}_{\bm{y}_0}^{\pi}(i)\|+\\
&~~~~~~~~~~~~~~~~~~~M_1 \|\bm{\phi}_{\bm{x}_0}^{\pi}(k)\|_{k>K}+M_1 \|\bm{\phi}_{\bm{y}_0}^{\pi}(k)\|_{k>K}\big)\\
&\leq \frac{\epsilon}{3}+\frac{\epsilon}{3}+\frac{\epsilon}{3}\leq \epsilon.
\end{split}
\end{equation}
Therefore, $W(\bm{x})$ is continuous over $B(\bm{0},\frac{r}{M})$.

For $\bm{x}_0\in \mathcal{R}_0$, let $L$ be the Lipschitz constant of $\ln(g(\bm{x})+1)$ over $\bm{x}\in \overline{X}$. Since $\mathcal{R}_0$ is open and $\bm{f}$ is Lipschitz continuous over $\bm{x}\in \mathbb{R}^n$ uniformly over $\bm{d}\in D$, we have that for $\epsilon$ satisfying $0<\epsilon <LK\delta'$, there exists an open neighborhood $O$ in $\mathcal{R}_0$ of $\bm{x}_0$ and $K>0$ such that \[\bm{\phi}_{\bm{y}_0}^{\pi}(k)\in B(\bm{0},\frac{r}{M}), \forall \bm{y}_0\in O, \forall \pi\in \mathcal{D}, \forall k\geq K\] and 
\[\|\bm{\phi}_{\bm{x}_0}^{\pi}(k)-\bm{\phi}_{\bm{y}_0}^{\pi}(k)\|\leq \frac{\epsilon}{LK}, \forall k\in [0,K],\] which implies that \[\|\bm{\phi}_{\bm{x}_0}^{\pi}(K)-\bm{\phi}_{\bm{y}_0}^{\pi}(K)\|\leq \frac{\epsilon}{LK}<\delta'.\]
Therefore, similar to the deduction in \eqref{rr1}, we have
\begin{equation*}
\begin{split}
&|W(\bm{x}_0)-W(\bm{y}_0)|\leq 2\epsilon.
\end{split}
\end{equation*}
In conclusion, $W(\bm{x}_0)$ is continuous over $\mathcal{R}_0$.

Next, we prove the continuity of $W'(\bm{x}_0)$. It is obvious that
\begin{equation*}
\begin{split}
&|W'(\bm{x}_0)-W'(\bm{y}_0)|\\
&\leq \sup_{\pi\in \mathcal{D}}\sup_{k\in \mathbb{N}}|\min_{j\in \{1,\ldots,n_{X}\}} h_{j,k}(\bm{x}_0,\pi)-\min_{j\in \{1,\ldots,n_{X}\}} h_{j,k}(\bm{y}_0,\pi)|.
\end{split}
\end{equation*}
As $\bm{x}_0\in \mathcal{R}_0$, $\lim_{k\rightarrow \infty}\min_{j\in \{1,\ldots,n_{X}\}} h_{j,k}(\bm{x}_0,\pi)=0$. Observing that $h_{j,k}$ is Lipschitz continuous over $\mathcal{R}_0$ and there exists $\beta(k):\mathbb{N}\rightarrow [0,\infty)$, which is independent of $\pi$, such that $\|\bm{\phi}_{\bm{x}_0}^{\pi}(k)\|\leq \beta(k)$ for $k\in \mathbb{N}, \pi\in \mathcal{D}$ and $\bm{x}_0\in \mathcal{R}_0$, we can find a neighborhood $B(\bm{x}_0,\delta)$ and a function $\gamma(k): \mathbb{N}\rightarrow [0,\infty)$ with $\lim_{k\rightarrow \infty}\gamma(k)= 0$ such that $|\min_{j\in \{1,\ldots,n_{X}\}}h_{j,k}(\bm{y}_0,\pi)|\leq \gamma(k)$ holds for $\bm{y}_0\in B(\bm{x}_0,\delta)$. This implies that the supremum 
\begin{equation*}
\begin{split}
\sup_{\pi\in \mathcal{D}}\sup_{k\in \mathbb{N}}&|\min_{j\in \{1,\ldots,n_{X}\}} h_{j,k}(\bm{x}_0,\pi)-\min_{j\in \{1,\ldots,n_{X}\}} h_{j,k}(\bm{y}_0,\pi)|
\end{split}
\end{equation*}
 is attained on a finite interval $[0,K]\cap \mathbb{N}$. On a compact time interval, the map $\bm{x}\rightarrow \min_{j\in \{1,\ldots,n_{X}\}}h_{j,k}(\bm{x},\pi)$ is Lipschitz continuous over $\bm{x}\in \mathcal{R}_0$ uniformly over $\pi\in \mathcal{D}$ since $h^X_{j}(\bm{x})$ and $\bm{f}(\bm{x},\bm{d})$ are Lipschitz continuous over $\bm{x}\in \mathcal{R}_0$ uniformly over $\bm{d}\in D$, implying that
\begin{equation*}
\begin{split}
\lim_{\bm{y}_0\rightarrow \bm{x}_0}\sup_{\pi\in \mathcal{D}}\sup_{k\in \mathbb{N}}&|\min_{j\in \{1,\ldots,n_{X}\}} h_{j,k}(\bm{x}_0,\pi)-\\
&\min_{j\in \{1,\ldots,n_{X}\}} h_{j,k}(\bm{y}_0,\pi)|=0.
\end{split}
\end{equation*}
This shows the desired continuity. 
  
The second assertion that $V(\bm{x})=\infty$ if $\bm{x}\notin \mathcal{R}_0$, can be proved by following the proof when $\bm{x}\notin \mathcal{R}_0$ in (a).

(c). From (b) we have that $V(\bm{x})=\infty$ for $\bm{x}\in \mathbb{R}^n\setminus \mathcal{R}_0$. Therefore, $v(\bm{x})=1$ for $\bm{x}\in \mathbb{R}^n\setminus \mathcal{R}_0$ due to the fact that $v(\bm{x})=1-e^{-V(\bm{x})}$ over $\mathbb{R}^n$. Therefore, $v(\bm{x})$ is continuous over $\mathbb{R}^n\setminus \mathcal{R}_0$. 

Also since $V(\bm{x})$ is continuous over $\mathcal{R}_0$, we have that $v(\bm{x})$ is continuous over $\mathcal{R}_0$. 

We just  prove that if $\lim_{\bm{x}\rightarrow \bm{y}}v(\bm{x})=v(\bm{y})$ for $\bm{x}\in \mathcal{R}_0$ and $\bm{y}\in \mathbb{R}^n\setminus \mathcal{R}_0$. According to (b) we have $\lim_{\bm{x}\rightarrow \bm{y}}V(\bm{x})=\infty$ and consequently 
$\lim_{\bm{x}\rightarrow \bm{y}}v(\bm{x})=1=v(\bm{y})$.

Above all, we have that $v(\bm{x})$ is continuous over $\mathbb{R}^n$. \qed
\end{pf}

Theorem \ref{conti} indicates that the interior of the maximal robust region of attraction can be obtained by computing either the value function $V(\bm{x})$ in \eqref{V0} or the value function $v(\bm{x})$ in \eqref{v}. Below we show that they can be computed by solving modified Bellman equations. For this sake, we first show that $V(\bm{x})$ and $v(\bm{x})$ satisfy the dynamic programming principle. 
\begin{lemma}
\label{dy}
Under Assumption \ref{app}, the following assertions are satisfied:
\begin{enumerate}
\item[(a)] For $\bm{x}\in \mathbb{R}^n$ and $k\in \mathbb{N}$, we have: 
\begin{equation}
\label{V11}
\begin{split}
&V(\bm{x})=\sup_{\pi\in \mathcal{D}}\max\big\{\\
&\sum_{i=1}^{k}g_{i-1}(\bm{x},\pi)+V(\bm{\phi}_{\bm{x}}^{\pi}(k)),\\
&\sup_{i\in [0,k-1]\cap \mathbb{N}}\{\sum_{j_1=1}^ig_{j_1-1}(\bm{x},\pi)-\min_{j\in \{1,\ldots,n_{X}\}} h_{j,i}(\bm{x},\pi)\}\\
&\big\}.
\end{split}
\end{equation}
\item[(b)] For $\bm{x}\in \mathbb{R}^n$ and $k\in \mathbb{N}$, we have: 
\begin{equation}
\label{v2}
\begin{split}
&v(\bm{x})=\sup_{\pi\in \mathcal{D}}\max\big\{\\
&1-\frac{1-v(\bm{\phi}_{\bm{x}}^{\pi}(k))}{\prod_{i=1}^{k}e^{g_{i-1}(\bm{x},\pi)}},\sup_{i\in [0,k-1]\cap \mathbb{N}}\{1-e^{-\overline{V}}\}\\
&\big\},
\end{split}
\end{equation}
where $
\overline{V}=\sum_{j_1=1}^i g_{j_1-1}(\bm{x},\pi)-\min_{j\in \{1,\ldots,n_{X}\}}h_{j,i}(\bm{x},\pi).$
\end{enumerate}
\end{lemma}
\begin{pf}
(a). Let 
\begin{equation}
\label{V1}
\begin{split}
&W(\bm{x}_0,k)=\sup_{\pi\in \mathcal{D}}\max\big\{\\
&\sum_{i=1}^{k}g_{i-1}(\bm{x}_0,\pi)+V(\bm{\phi}_{\bm{x}_0}^{\pi}(k)),\\
&\sup_{i\in [0,k-1]\cap \mathbb{N}} \{\sum_{j_1=1}^ig_{j_1-1}(\bm{x}_0,\pi)-\min_{j\in \{1,\ldots,n_{X}\}} h_{j,i}(\bm{x}_0,\pi)\}\\
&\big\}.
\end{split}
\end{equation}
We will prove that $|W(\bm{x}_0,k)-V(\bm{x}_0)|\leq \epsilon$, $\forall \epsilon>0$.

From \eqref{V0}, for any $\epsilon_1>0$, there exists $\pi\in \mathcal{D}$ such that 
\begin{equation*}
\begin{split}
&V(\bm{x}_0)\leq \epsilon_1+
\\
&~~~~~~~\sup_{k\in \mathbb{N}}\big\{\sum_{i=1}^k g_{i-1}(\bm{x}_0,\pi)-\min_{j\in \{1,\ldots,n_{X}\}} h_{j,k}(\bm{x}_0,\pi)\big\}\}.
\end{split}
\end{equation*}

We respectively define $\pi_1\in \mathcal{D}$ and $\pi_2 \in \mathcal{D}$ as follows:  $\pi_1(i)=\pi(i)$ for $i=0,\ldots,k-1$, and $\pi_2(i)=\pi(i+k)$ for $i\in \mathbb{N}$, and $\bm{y}=\bm{\phi}_{\bm{x}_0}^{\pi_1}(k)$, then obtain that 
\begin{equation*}
\begin{split}
&W(\bm{x}_0,k)\geq \max\big\{\\
&\sum_{i=1}^{k}g_{i-1}(\bm{x}_0,\pi)+V(\bm{y}),\\
&\sup_{i\in [0,k-1]\cap\mathbb{N}}\{\sum_{j_1=1}^ig_{j_1-1}(\bm{x}_0,\pi)-\min_{j\in \{1,\ldots,n_{X}\}}h_{j,i}(\bm{x}_0,\pi)\}\\
&\big\}\\
&\geq\max\big\{\\
&\sum_{i=1}^{k}g_{i-1}(\bm{x}_0,\pi_{1})+\\
&\sup_{l\in [k,\infty)\cap \mathbb{N}}\{ \sum_{i'=1}^{l-k}g_{i'-1}(\bm{y},\pi_{2})-\min_{j\in \{1,\ldots,n_{X}\}}h_{j,l-k}(\bm{y},\pi_{2})\},\\
&\sup_{i\in [0,k-1]\cap \mathbb{N}} \{\sum_{j_1=1}^ig_{j_1-1}(\bm{x}_0,\pi_{1})-\min_{j\in \{1,\ldots,n_{X}\}} h_{j,i}(\bm{x}_0,\pi_{1})\}\\
&\big\}\\
&\geq\max\big\{\\
&\sup_{l\in [k,\infty)\cap \mathbb{N}} \{\sum_{j_1=1}^{l}g_{j_1-1}(\bm{x}_0,\pi)-\min_{j\in \{1,\ldots,n_{X}\}}h_{j,l}(\bm{x}_0,\pi)\}, \\
&\sup_{i\in [0,k-1]\cap \mathbb{N}} \{\sum_{j_1=1}^ig_{j_1-1}(\bm{x}_0,\pi)-\min_{j\in \{1,\ldots,n_{X}\}}h_{j,i}(\bm{x}_0,\pi)\}\\
&\big\}\\
&\geq V(\bm{x}_0)-\epsilon_1.
\end{split}
\end{equation*}
Therefore,  $V(\bm{x}_0)\leq W(\bm{x}_0,k)+\epsilon_1.$

According to \eqref{V1}, for any $\epsilon_1>0$, there exists a perturbation input policy $\pi_{1}\in \mathcal{D}$ such that 
\begin{equation*}
\begin{split}
&W(\bm{x}_0,k)\leq \epsilon_1+\max\big\{\\
&\sum_{i=1}^{k}g_{i-1}(\bm{x}_0,\pi_{1})+V(\bm{\phi}_{\bm{x}_0}^{\pi_{1}}(k)),\\
&\sup_{i\in [0,k-1]\cap \mathbb{N}} \{\sum_{j_1=1}^ig_{j_1-1}(\bm{x}_0,\pi_{1})-\min_{j\in \{1,\ldots,n_{X}\}}h_{j,i}(\bm{x}_0,\pi_{1})\}\\
&\big\}.
\end{split}
\end{equation*}
Also, by the definition of $V$, i.e. \eqref{V0}, for any $\epsilon_1$, there exists an input policy $\pi_{2}\in \mathcal{D}$ such that 
\begin{equation*}
\begin{split}
V(\bm{y})&\leq \epsilon_1+\\
&\sup_{l\in \mathbb{N}}\big\{\sum_{i=1}^l g_{i-1}(\bm{y},\pi_{2})-\min_{j\in \{1,\ldots,n_{X}\}}h_{j,l}(\bm{y},\pi_{2})\big\},
\end{split}
\end{equation*}
where $\bm{y}=\bm{\phi}_{\bm{x}_0}^{\pi_{1}}(k)$. We define $\pi$:
\begin{equation*}
   \pi(i) = \left\{
   \begin{array}{ll}
     \pi_1(i), & i\in [0,k)\cap \mathbb{N} \\
     \pi_2(i-k), & i\in [k,\infty)\cap \mathbb{N}
   \end{array}
   \right..
\end{equation*}
Therefore, we infer that 
\begin{equation*}
\begin{split}
&W(\bm{x}_0,k)\leq \epsilon_1+\max\big\{\\
&\sum_{i=1}^{k}g_{i-1}(\bm{x}_0,\pi_{1})+V(\bm{y}),\\
&\sup_{i\in [0,k-1]\cap \mathbb{N}}\{\sum_{j_1=1}^ig_{j_1-1}(\bm{x}_0,\pi_{1})-\min_{j\in \{1,\ldots,n_{X}\}}h_{j,i}(\bm{x}_0,\pi_{1})\}\\
&\big\}\\
&\leq 2\epsilon_1+\max\big\{\\
&\sum_{i=1}^{k}g_{i-1}(\bm{x}_0,\pi_{1})+\\
&\sup_{l\in [k,\infty)\cap \mathbb{N}}\{ \sum_{i'=1}^{l-k}g_{i'-1}(\bm{y},\pi_{2})-\min_{j\in \{1,\ldots,n_{X}\}} h_{j,l-k}(\bm{y},\pi_{2}), \\
&\sup_{i\in [0,k-1]\cap \mathbb{N}} \{\sum_{j_1=1}^ig_{j_1-1}(\bm{x}_0,\pi_{1})-\min_{j\in \{1,\ldots,n_{X}\}} h_{j,i}(\bm{x}_0,\pi_{1})\}\\
&\big\}\\
&\leq 2\epsilon_1+\max\big\{\\
&\sup_{l\in [k,\infty)\cap \mathbb{N}}\{\sum_{j_1=1}^l g_{j_1-1}(\bm{x}_0,\pi)-\min_{j\in \{1,\ldots,n_{X}\}} h_{j,l}(\bm{x}_0,\pi)\},\\
&\sup_{i\in [0,k-1]\cap \mathbb{N}}\{\sum_{j_1=1}^ig_{j_1-1}(\bm{x}_0,\pi)-\min_{j\in \{1,\ldots,n_{X}\}}h_{j,i}(\bm{x}_0,\pi)\} \\
&\big\}\\
&\leq V(\bm{x}_0)+2\epsilon_1.
\end{split}
\end{equation*}

Therefore, we finally have $|W-V|\leq \epsilon=2\epsilon_1$, implying that $V=W$ since $\epsilon_1$ is arbitrary.

(b). \eqref{v2} can be obtained using $v(\bm{x}_0)=1-e^{-V(\bm{x}_0)}$. \qed
\end{pf}

Based on Lemma \ref{dy} we can infer that the value functions $V(\bm{x}_0)$ and $v(\bm{x}_0)$ are solutions to the two modified Bellman equations \eqref{eq1} and \eqref{e-v}, respectively. 
\begin{theorem}
\label{equations}
Under Assumption \ref{app}, the value function $V$ is the unique continuous solution to the modified Bellman equation 
\begin{equation}
\label{eq1}
\begin{split}
&\min\big\{\inf_{\bm{d}\in D}\{V-V(\bm{f})-\ln(g+1)\},\\
&~~~~~V+\min_{j\in\{1,\ldots,n_{X}\}}\ln(l(1-h^X_j))\big\}=0, \forall \bm{x}\in \mathcal{R}_0,\\
&V(\bm{0})=0.
\end{split}
\end{equation}
The value function $v$ is the unique bounded and continuous solution to the modified Bellman equation
\begin{equation}
\label{e-v}
\begin{split}
&\min\big\{\inf_{\bm{d}\in D}\{v-v(\bm{f})-g\cdot(1-v)\},\\
&~~~~~v-1+\min_{j\in \{1,\ldots,n_{X}\}} l(1-h^X_j)\big\}=0, \forall \bm{x}\in \mathbb{R}^n,\\
&v(\bm{0})=0.\\
\end{split}
\end{equation}
\end{theorem}
\begin{pf}
The fact that the value functions $V(\bm{x})$ in \eqref{V0} and $v(\bm{x})$ in \eqref{v} are solutions to \eqref{eq1} and \eqref{e-v} respectively  can be verified when $k=1$ in \eqref{V11} and \eqref{v2}.  

Here, we just prove the uniqueness of solutions to \eqref{e-v}. The uniqueness of solution to \eqref{eq1} can be guaranteed by the  relationship $v(\bm{x})=1-e^{-V(\bm{x})}$ for $\bm{x}\in \mathbb{R}^n$. 

Assume that $\tilde{v}$ is a bounded and continuous solution to \eqref{e-v} as well, we need to prove that $v=\tilde{v}$ over $\bm{x}\in \mathbb{R}^n$, where $v<1$ over $\mathcal{R}_0$ and $v=1$ over $\mathbb{R}^n\setminus \mathcal{R}_0$. Assume that there exists $\bm{y}_0$ such that $\tilde{v}(\bm{y}_0)\neq v(\bm{y}_0)$. First let's assume $v(\bm{y}_0)>\tilde{v}(\bm{y}_0)$ and $v(\bm{y}_0)\geq 1$. Obviously, $\bm{y}_0\neq \bm{0}$ and consequently $g(\bm{y}_0)>0$. Since both $v$ and $\tilde{v}$ satisfy \eqref{e-v}, we have that $$\inf_{\bm{d}\in D}\{v(\bm{y}_0)-v(\bm{f}(\bm{y}_0,\bm{d}))-g(\bm{y}_0)(1-v(\bm{y}_0))\}=0.$$ Since $v$ is continuous over $\mathbb{R}^n$ and $\bm{f}$ is continuous over $\mathbb{R}^n\times D$, there exists $\bm{d}'_1\in D$ such that $v(\bm{y}_0)-v(\bm{f}(\bm{y}_0,\bm{d}'_1))-g(\bm{y}_0)(1-v(\bm{y}_0))=0$. Since $\tilde{v}(\bm{y}_0)-\tilde{v}(\bm{f}(\bm{y}_0,\bm{d}'_1))-g(\bm{y}_0)(1-\tilde{v}(\bm{y}_0))\geq 0$, we obtain that 
\begin{equation*}
\begin{split}
v(\bm{f}(\bm{y}_0,\bm{d}'_1))-&\tilde{v}(\bm{f}(\bm{y}_0,\bm{d}'_1))\\
&\geq (v(\bm{y}_0)-\tilde{v}(\bm{y}_0))(1+g(\bm{y}_0)).
\end{split}
\end{equation*}
Let $\bm{y}_1=\bm{\phi}_{\bm{y}_0}^{\pi_{1}}(1)$, where $\pi_1(0)=\bm{d}'_1$, then $v(\bm{y}_1)>\tilde{v}(\bm{y}_1)$. Also, we have $v(\bm{y}_0)\leq v(\bm{y}_1)$. Moreover, $\bm{y}_1 \neq \bm{0}$, $g(\bm{y}_1)>0$. We continue the above deduction from $\bm{y}_0$ to $\bm{y}_1$, and obtain that there exists $\bm{d}'_2 \in D$ such that 
\begin{equation*}
\begin{split}
v(\bm{f}(\bm{y}_1,\bm{d}'_2))-&\tilde{v}(\bm{f}(\bm{y}_1,\bm{d}'_2))\\
&\geq (v(\bm{y}_1)-\tilde{v}(\bm{y}_1))(1+g(\bm{y}_1)).
\end{split}
\end{equation*}
Thus, we have 
\begin{equation*}
\begin{split}
&v(\bm{f}(\bm{y}_1,\bm{d}'_2))-\tilde{v}(\bm{f}(\bm{y}_1,\bm{d}'_2))\geq \\
&(v(\bm{y}_0)-\tilde{v}(\bm{y}_0))(1+g(\bm{y}_1))(1+g(\bm{y}_0)).
\end{split}
\end{equation*}
Let $\bm{y}_2=\bm{\phi}_{\bm{y}_1}^{\pi_{2}}(1)$, where $\pi_2(0)=\bm{d}'_2$, then $v(\bm{y}_2)>\tilde{v}(\bm{y}_2)$. Also, $v(\bm{y}_1)\leq v(\bm{y}_2)$.

Analogously, we deduce that for $k\in \mathbb{N}$,
\begin{equation*}
\begin{split}
v(\bm{f}&(\bm{y}_k,\bm{d}'_{k+1}))-\tilde{v}(\bm{f}(\bm{y}_k,\bm{d}'_{k+1}))\geq\\
& (v(\bm{y}_0)-\tilde{v}(\bm{y}_0))(1+g(\bm{y}_k))\cdots(1+g(\bm{y}_0)).
\end{split}
\end{equation*}
Moreover, let $\bm{y}_{k+1}=\bm{\phi}_{\bm{y}_k}^{\pi_{{k+1}}}(1)$, then $v(\bm{y}_k)\leq v(\bm{y}_{k+1})$, where $\pi_{k+1}(0)=\bm{d}'_{k+1}$. This implies that $\lim_{k\rightarrow \infty} \bm{y}_k\neq \bm{0}$ and thus $\bm{y}_k\notin B(\bm{0},\overline{\epsilon})$ for $k\in \mathbb{N}$, where $B(\bm{0},\overline{\epsilon})$ is defined in \eqref{half}.
Assume that $c_0=\inf\{g(\bm{x})\mid \bm{x}\in\mathbb{R}^n\setminus B(\bm{0},\overline{\epsilon})\}$. Obviously, $c_0>0$. Therefore, 
\begin{equation*}
\begin{split}
v(\bm{f}(\bm{y}_k,\bm{d}'_{k+1}))-&\tilde{v}(\bm{f}(\bm{y}_k,\bm{d}'_{k+1}))\\
&\geq (v(\bm{y}_0)-\tilde{v}(\bm{y}_0))(1+c_0)^{k+1},
\end{split}
\end{equation*}
 implying that $\lim_{k\rightarrow \infty}v(\bm{y}_k)=\infty$, which contradicts the fact that $v$ is bounded over $\mathbb{R}^n$. 

Next, assume $v(\bm{y}_0)>\tilde{v}(\bm{y}_0)$ and $v(\bm{y}_0)<1$. According to Theorem \ref{conti}, every possible trajectory starting from $\bm{y}_0$ will eventually approach $\bm{0}$. Also, we have $$\inf_{\bm{d}\in D}\{v(\bm{y}_0,\bm{d})-v(\bm{f}(\bm{y}_0,\bm{d}))-g(\bm{y}_0)(1-v(\bm{y}_0))\}=0.$$ Following the deduction mentioned above, we have $$v(\bm{y}_k)-\tilde{v}(\bm{y}_{k})\geq v(\bm{y}_0)-\tilde{v}(\bm{y}_0),\forall k\in \mathbb{N}.$$ Since $\lim_{k\rightarrow \infty}\tilde{v}(\bm{y}_{k})=0$, $\lim_{k\rightarrow \infty}v(\bm{y}_k)\geq v(\bm{y}_0)-\tilde{v}(\bm{y}_0)$ holds, contradicting $\lim_{k\rightarrow \infty}v(\bm{y}_k)=0$. 

For the case that $\tilde{v}(\bm{y}_0)>v(\bm{y}_0)$, we can obtain similar contradictions by following the proof procedure mentioned above with $v$ and $\tilde{v}$ reversed. \qed
\end{pf}

\end{document}